
\input epsf.tex
\input amssym.def
\input amssym
\magnification=1100
\baselineskip = 0.22truein
\lineskiplimit = 0.01truein
\lineskip = 0.01truein
\vsize = 8.5truein
\voffset = 0.2truein
\parskip = 0.10truein
\parindent = 0.3truein
\settabs 12 \columns
\hsize = 5.8truein
\hoffset = 0.2truein

\setbox\strutbox=\hbox{%
\vrule height .708\baselineskip
depth .292\baselineskip
width 0pt}
\font\caps=cmcsc10
\font\bigtenrm=cmr10 at 14pt

\def\sqr#1#2{{\vcenter{\vbox{\hrule height.#2pt
\hbox{\vrule width.#2pt height#1pt \kern#1pt
\vrule width.#2pt}
\hrule height.#2pt}}}}
\def\square{\mathchoice\sqr46\sqr46\sqr{3.1}6\sqr{2.3}4}

\centerline{\bigtenrm THE CROSSING NUMBER OF COMPOSITE KNOTS}
\tenrm
\vskip 14pt
\centerline{MARC LACKENBY}
\vskip 18pt
\centerline{\caps 1. Introduction}
\vskip 6pt

One of the most basic questions in knot theory remains
unresolved: is crossing number additive under connected
sum? In other words, does the equality $c(K_1 \sharp K_2) =
c(K_1) + c(K_2)$ always hold, where $c( K )$ denotes the
crossing number of a knot $K$ and $K_1 \sharp K_2$ is the connected
sum of two (oriented) knots $K_1$ and $K_2$? The inequality
$c(K_1 \sharp K_2) \leq c(K_1) + c(K_2)$ is trivial, but very little
more is known in general. Equality has been
established for certain classes of knots, most notably
when $K_1$ and $K_2$ are both alternating ([3], [6], [7])
and when $K_1$ and $K_2$ are both torus knots [1].
In this paper, we provide the first non-trivial lower
bound on $c(K_1 \sharp K_2)$ that applies to all knots
$K_1$ and $K_2$.

\noindent {\bf Theorem 1.1.} {\sl Let $K_1, \dots, K_n$ be oriented
knots in the 3-sphere. Then
$${c(K_1) + \dots + c(K_n) \over 152} \leq c(K_1 \sharp \dots \sharp K_n)
\leq c(K_1) + \dots + c(K_n).$$}

More generally, one can speculate about the crossing
number of satellite knots. Here, there are a variety
of conjectures, all of which remain wide open at present.
The simplest of these asserts that the crossing number
of a non-trivial satellite knot is at least the crossing
number of its companion. To explain this, we fix some
terminology. A knot $K$ is a {\sl non-trivial satellite
knot with companion knot $L$} if $K$ lies in a regular neighbourhood
$N(L)$ of the non-trivial knot $L$, and $K$ does not lie in a 3-ball
contained in $N(L)$, and $K$ is not a core curve of the solid
torus $N(L)$. In a forthcoming article [4], we will prove the following
result, by generalising the methods in this paper.

\noindent {\bf Theorem 1.2.} {\sl There is a universal computable constant
$N \geq 1$ with the following property. Let $K$ be a non-trivial
satellite knot, with companion knot $L$. Then $c(K) \geq c(L)/N$. }

Here, `universal' means that $N$ is just a number, and
`computable' means that we have an algorithm to determine it.
However, the constant $N$ is more difficult to calculate than in the case
of composite knots. We hope to find an explicit upper bound on $N$,
but it will probably be significantly bigger than 152.

\vfill\eject
Let us start with an
outline of the proof of Theorem 1.1. Let $K_1, \dots, K_n$ be
a collection of oriented knots. Our aim is to
show that $c(K_1) + \dots + c(K_n) \leq 152 \, c(K_1 \sharp \dots
\sharp K_n)$. It is not hard to show
that we may assume that each $K_i$ is prime and non-trivial.
Let $D$ be a diagram of $K_1 \sharp \dots
\sharp K_n$ having minimal crossing number.
Our goal is to construct a diagram $D'$ for
the distant union $K_1 \sqcup \dots \sqcup K_n$
such that $c(D') \leq 152 \, c(D)$. (The {\sl distant union}
of oriented knots $K_1, \dots, K_n$, denoted $K_1 \sqcup \dots \sqcup K_n$,
is constructed by starting with $n$ disjoint 3-balls in the 3-sphere,
and for $i = 1, \dots, n$, placing a copy of $K_i$ in the $i$th ball.)
Theorem 1.1 is then a consequence of the following easy lemma
which is proved in Section 2.

\noindent {\bf Lemma 2.1.} {\sl Let $K_1 \sqcup \dots \sqcup K_n$
be the distant union of oriented knots $K_1, \dots, K_n$.
Then
$$c(K_1 \sqcup \dots \sqcup K_n) = c(K_1) + \dots + c(K_n).$$}

So, the key to the proof of Theorem 1.1 is
to construct the diagram $D'$ with
$c(D') \leq 152 \, c(D)$. 
Let $X$ be the exterior of $K = K_1 \sharp \dots \sharp K_n$.
Arising from the connected sum construction of $K$, there
is a collection of $n$ disjoint annuli $A_1, \dots, A_n$ properly
embedded in $X$. These are shown in Figure 1. Let $A$ be $A_1 \cup \dots \cup A_n$.

\vskip 18pt
\centerline{
\epsfxsize=3in
\epsfbox{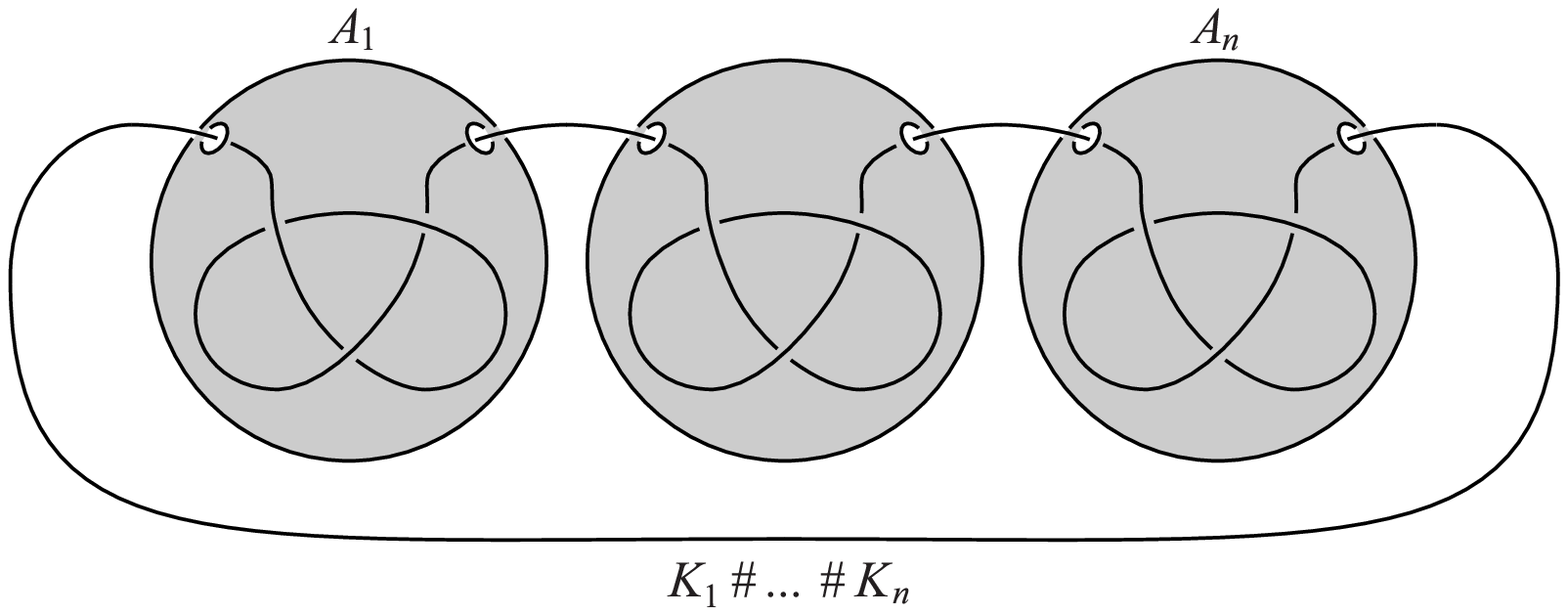}
}
\vskip 6pt
\centerline{Figure 1.}

If one were to cut $X$ along $A$, the resulting
3-manifold would be the disjoint union of 3-manifolds $X_1, \dots, X_n$ and $Y$,
where each $X_i$ is homeomorphic to the exterior of $K_i$,
and $Y$ is the component with a copy of each of $A_1, \dots, A_n$ in
its boundary. Each $X_i$ is separated from the other components by a 2-sphere,
which is made up of $A_i$ and two meridian discs for $K$. 
Thus, if one were to choose on the boundary of each $X_i$ a simple closed curve $C_i$
that intersects the meridian of $K_i$ just once, then the union of
these simple closed curves would be $K_1 \sqcup \dots \sqcup K_n$.

Now, the given diagram $D$ for $K_1 \sharp \dots \sharp K_n$ need not look
like that shown in Figure 1. So, the annuli $A$ may sit inside $S^3$
in a complicated way. Hence, the 3-manifolds $X_1, \dots, X_n$
may be embedded in $S^3$ in a highly twisted fashion, which
means that {\sl a priori} the diagram $D'$ of $K_1 \sqcup \dots \sqcup K_n$
obtained by projecting $C_1 \cup \dots \cup C_n$ may be very complex. In particular, it may have many more crossings
than $D$. The goal is to gain enough
control over the annuli $A$ and hence over the manifolds $X_1, \dots, X_n$,
so that we can bound the number of crossings in $D'$. For this,
our main tool will be normal surface theory.

We first construct a handle structure for $X$, arising from
the diagram $D$. Since the annuli are essential in $X$,
they may be placed in normal form with respect
to this handle structure. This alone does not give
us enough control, since the annuli may run through
each handle many times. However, if a handle
contains many normal discs, they fall into a bounded
number of disc types so that any
two discs of the same type are normally parallel. Our
main technical achievement is to show that
the curves $C_i$ can be chosen so that they
miss the normal discs that have parallel copies on both
sides. So, they run through each handle a bounded number of
times and in a controlled way.
Thus, we can bound the number of new crossings
that are introduced when constructing $D'$ from $D$,
to obtain the inequality $c(D') \leq 152 \, c(D)$.
The constant 152 arises from the combinatorics of
normal discs in our chosen handle structure.

The paper is organised as follows. In Section 2, we will prove Lemma 2.1,
which gives the formula for $c(K_1 \sqcup \dots \sqcup K_n)$. In Section 3,
we assign a handle structure to $X$
using the diagram $D$ for $K$. In Section 4, we recall some of the theory of
normal surfaces in handle structures. (For a more complete reference, see [5].)
Once we have placed the annuli $A$ into normal form in $X$, we cut
$X$ along $A$ and discard the component $Y$ that contains a copy of $A$ in
its boundary. The resulting 3-manifold $M$ is a disjoint union of
$X_1, \dots, X_n$. It inherits a handle structure. Let $S$ be the
copy of $A$ in $\partial M$. In Section 5, we define
the notion of a generalised parallelity bundle in $M$. This is a subset
${\cal B}$ of $M$ homeomorphic to an $I$-bundle over a surface, such that the
$\partial I$-bundle is ${\cal B} \cap S$, and which has various other properties.
An example of a generalised parallelity bundle is the union of all
the handles of $M$ that lie between parallel normal discs of $A$.
In Section 5, we establish the existence of a generalised parallelity bundle
${\cal B}$ that contains all these handles, and maybe others, and which has the following key property:
{\sl either} the handle structure of $M$ admits a certain type of simplification,
known as an annular simplification (in which case, we perform this simplification
and continue) {\sl or} each component of ${\cal B}$
is an $I$-bundle over a disc or has incompressible vertical boundary.
(The vertical boundary of an $I$-bundle is the
closure of the subset of the boundary that does not lie in the $\partial I$-bundle;
it is a collection of annuli.) In Section 6, we complete the proof of Theorem 1.1.
We note that $M$ does not admit any essential embedded
annuli with boundary lying in $S$, because we may assume that the knots $K_1, \dots, K_n$ are prime.
Thus, with further work, we deduce that the generalised parallelity bundle ${\cal B}$
is a collection of $I$-bundles over discs. Hence, its $\partial I$-bundle
does not separate the boundary components of each component of $S$. It is therefore possible
to choose the curves $C_1, \dots, C_n$ so that they avoid ${\cal B}$. 
In particular, they avoid the handles of $M$ that lie between parallel normal
discs of $A$. They therefore they run through each handle 
of $X$ a bounded number of times and in a controlled way. The final parts
of Section 6 are devoted to quantifying this control, and justifying the
constant $152$.

\vskip 18pt
\centerline{\caps 2. The crossing number of the distant
union of knots}
\vskip 6pt

In this short section, we prove the following
lemma, which is a key step in the proof of Theorem 1.1.

\noindent {\bf Lemma 2.1.} {\sl Let $K_1 \sqcup \dots \sqcup K_n$
be the distant union of oriented knots $K_1, \dots, K_n$.
Then
$$c(K_1 \sqcup \dots \sqcup K_n) = c(K_1) + \dots + c(K_n).$$}

\noindent {\sl Proof.} To prove the inequality
$c(K_1 \sqcup \dots \sqcup K_n) \geq c(K_1) + \dots + c(K_n)$,
consider a diagram $D$ of 
$K_1 \sqcup \dots \sqcup K_n$ with minimal crossing number.
From this, one can construct a diagram $D_i$ of
$K_i$, by eliminating all other components.
Thus, $c(D_i)$ is the number of crossings of
$D$ where the over-arc and under-arc both lie
in $K_i$. The sum $\sum_{i=1}^n c(D_i)$ therefore enumerates
a subset of the crossings of $D$. Hence,
$$c(K_1 \sqcup \dots \sqcup K_n) = c(D) 
\geq \sum_{i=1}^n c(D_i) \geq \sum_{i=1}^n c(K_i).$$ 
The inequality in the other direction is trivial, since one can construct
a diagram for $K_1 \sqcup \dots \sqcup K_n$ from
minimal crossing number diagrams of $K_1, \dots, K_n$.
$\square$

It is intriguing that distant unions are so much 
more tractable than connected sums. 

\vskip 18pt
\centerline{\caps 3. A handle structure from a diagram}
\vskip 6pt

In this section, we describe a method for constructing
a handle structure ${\cal H}$ on the exterior of a knot $K$, starting
with a diagram $D$ for $K$. In some sense, it not particularly
important how one does this. As long as one picks a handle structure
in a reasonably sensible way, then the remaining techniques in
this paper will give a result like Theorem 1.1, but possibly with
152 replaced by a different constant.

The diagram $D$ is
a 4-valent graph embedded in a 2-sphere $S^2$, 
with crossing information at each vertex.
Associated with $D$, there are two collections of disjoint
arcs embedded in the 2-sphere, which we denote by $D_+$
and $D_-$. Roughly speaking, $D_+$ is the collection of
arcs made by the pen when one draws the knot. That is,
one makes two small cuts near each vertex of $D$, so that
the over-arc runs smoothly through the crossing, but the
under-arcs are terminated. The resulting collection of arcs is $D_+$.
The arcs $D_-$ are defined similarly, but where the over-arcs are cut at each
crossing and the under-arcs run through smoothly.

We realise the 3-sphere as the set of points $(x_1, x_2, x_3, x_4)$
in ${\Bbb R}^4$ with Euclidean norm 2, say. We embed
the diagram 2-sphere as the equator $\{ x_4 = 0 \}$.
The north and south poles of $S^3$ are the points
$(0,0,0,2)$ and $(0,0,0,-2)$ respectively.
There is a homeomorphism from the complement of
these two points to $S^2 \times (-2,2)$,
such that projection $S^2 \times (-2,2) \rightarrow
(-2,2)$ onto the second factor of the product
agrees with the height function $x_4$.

The diagram $D$ specifies an embedding of
$K$ into the 3-sphere, as follows. Away from
a small regular neighbourhood of the crossings,
the knot lies in the diagram 2-sphere.
Near each crossing, the knot leaves this
2-sphere, forming two arcs, one lying
above the diagram, and one below it.
Specifically, the over-arc runs
vertically up from the diagram, then
runs horizontally at height $x_4 = 1$
say, and then goes vertically back down
to the diagram. The under-arc has a similar itinerary below the diagram.
Thus, the diagrammatic projection map from the complement of the
north and south poles onto the diagram 2-sphere
is the product projection map $S^2 \times (-2,2) \rightarrow S^2$.

We pick a point $\infty$ in the diagram 2-sphere $S^2$ that is
distant from the crossings, and assign a Euclidean metric
to $S^2 - \{ \infty \}$.

We now define a handle structure ${\cal H}'$ on the exterior of $K$.
The handle structure ${\cal H}$ that we actually use in the proof of Theorem 1.1
will be a slight modification of this.

We start with the 0-handles of ${\cal H}'$. Near each crossing,
we place four 0-handles, as shown in Figure 2.
Instead of using round 3-balls for the 0-handles, it
is slightly more convenient to take each to be of the form
$D^2 \times [-1, 1]$, where $D^2$ is a Euclidean disc, and the second factor
is the $x_4$ co-ordinate.

We now add the 1-handles. Near each crossing, 
we add four 1-handles. These are `horizontal',
in the sense that they are regular neighbourhoods
of arcs in the diagram 2-sphere. These four 1-handles
run between the four 0-handles like
the edges of a square. Note that these
1-handles do indeed lie in the exterior of
$K$ because $K$ skirts above and below the
diagram 2-sphere at these points. In addition, for each
edge of the 4-valent graph of the knot projection,
we add two horizontal 1-handles, which lie either side of
the edge and run parallel to it.

\vskip 18pt
\centerline{
\epsfxsize=3.7in
\epsfbox{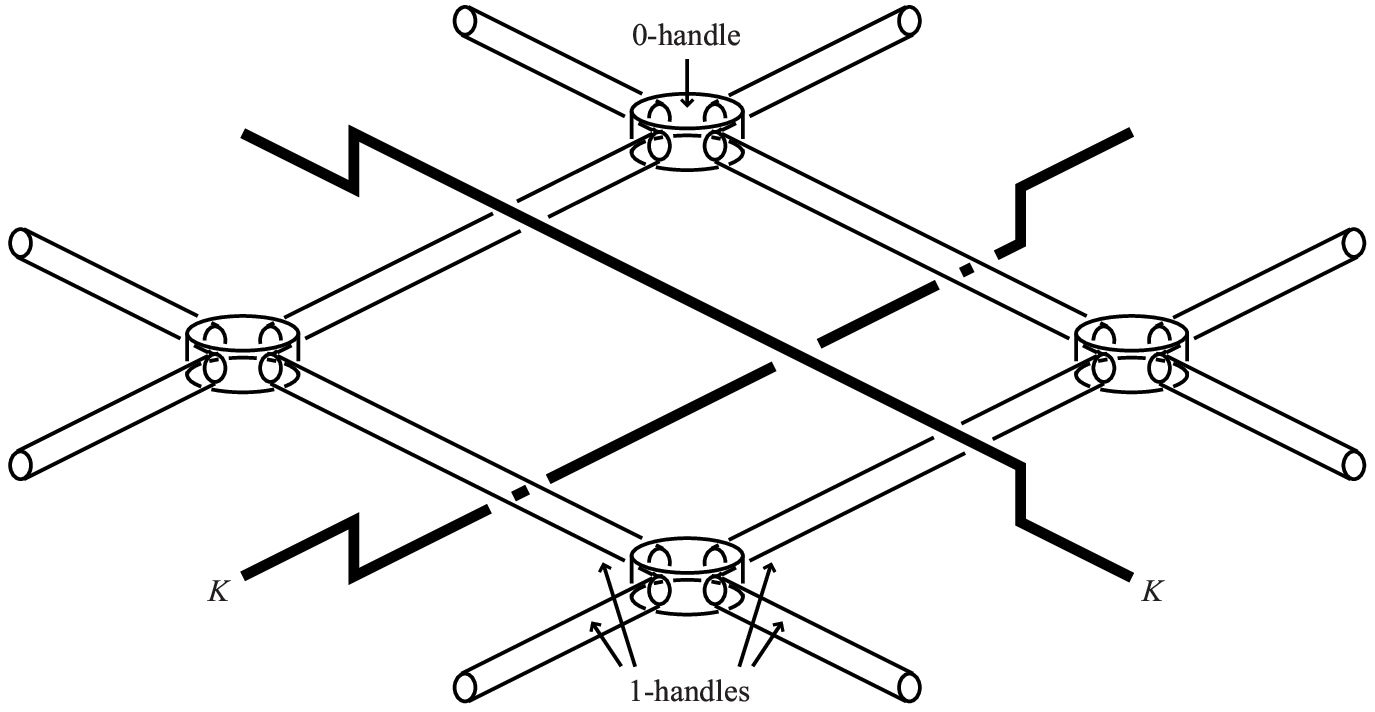}
}
\vskip 6pt
\centerline{Figure 2.}

We now specify where the 2-handles lie.
A square-shaped 2-handle, as shown in Figure 3, is attached
to the square-shaped configuration of 1-handles and 0-handles
near each crossing. It is 
`horizontal', in the sense that it is
a thin regular neighbourhood of a subset of
the diagram 2-sphere. Thus, it is attached to the
1-handles and 0-handles in the `plane of the diagram'.

\vskip 18pt
\centerline{
\epsfxsize=3.7in
\epsfbox{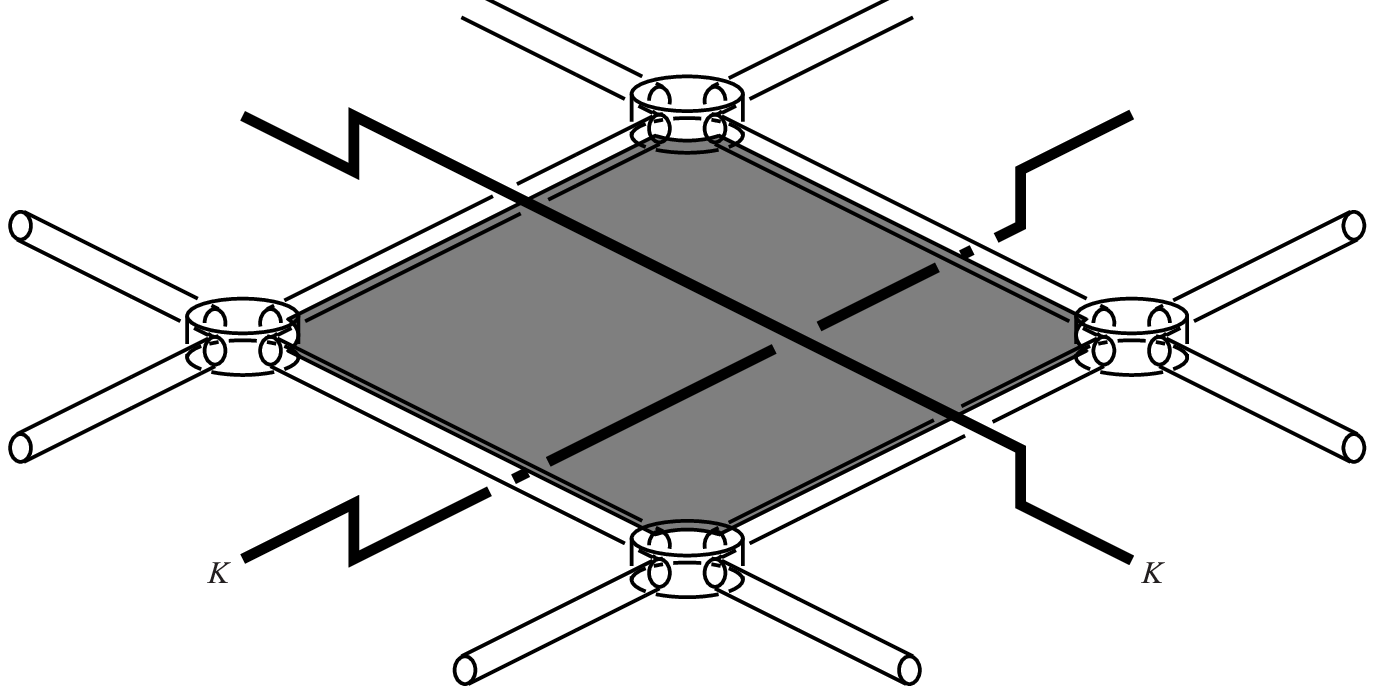}
}
\vskip 6pt
\centerline{Figure 3.}

Associated with each region of the diagram, there
is also a horizontal 2-handle. Its attaching annulus runs
along the 1-handles and 0-handles that lie within that
region, as in Figure 4. It too is attached to the
1-handles and 0-handles in the `plane of the diagram'.

\vskip 18pt
\centerline{
\epsfxsize=3.7in
\epsfbox{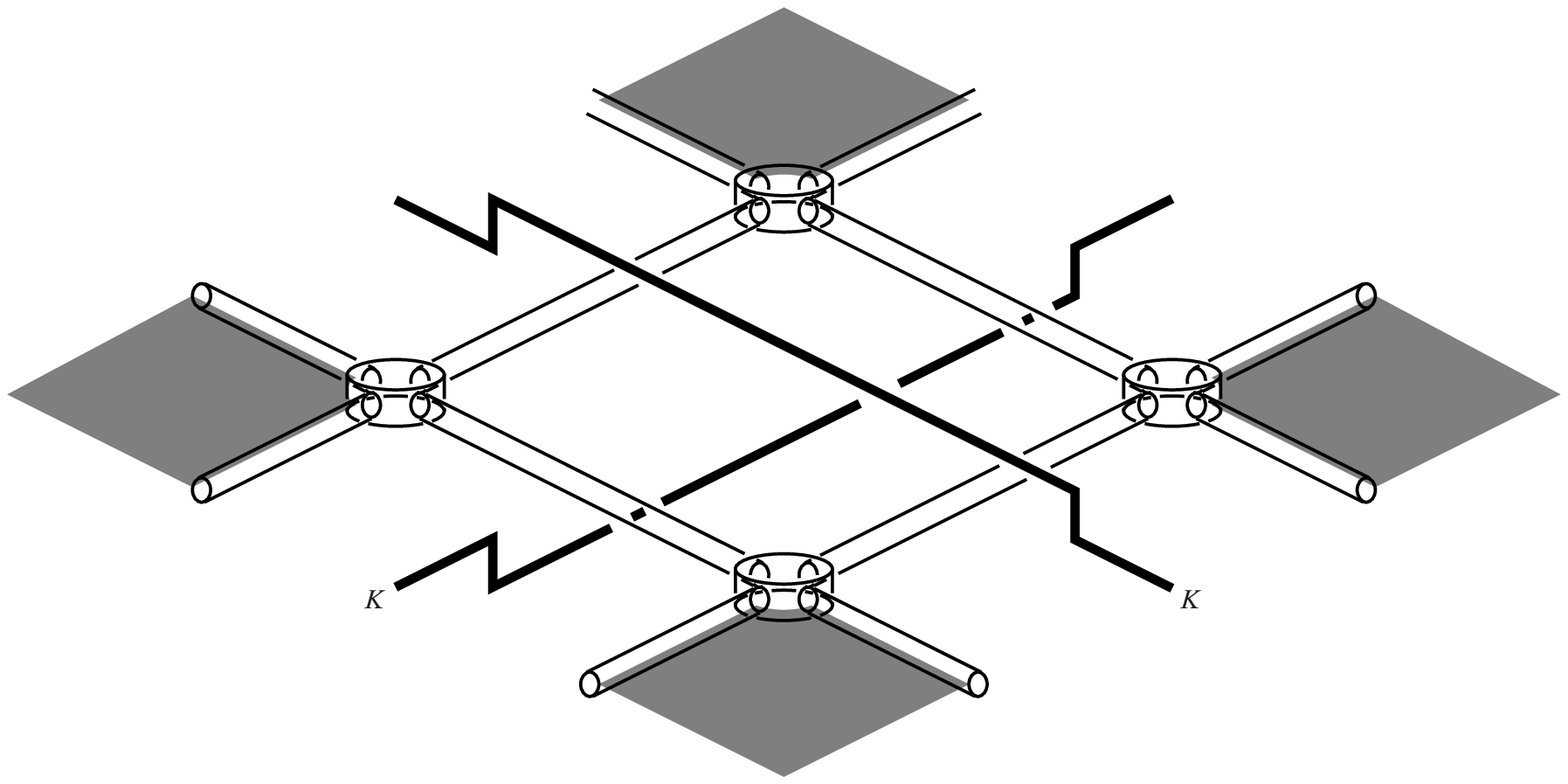}
}
\vskip 6pt
\centerline{Figure 4.}

There is one final type of 2-handle. Each is associated with
an arc of $D_+$ or $D_-$. It is attached along the 1-handles
and 0-handles that encircle this
arc, and the 2-handle itself lies over the diagram (in the case of
$D_+$) or under the diagram (in the case of $D_-$). Part of
such a 2-handle is shown in Figure 5. We may suppose that
the points where the 2-handle is attached to the 0-handles and 1-handles
lie just above the diagram (in the case of $D_+$) or just below
the diagram (in the case of $D_-$).

\vskip 18pt
\centerline{
\epsfxsize=3.7in
\epsfbox{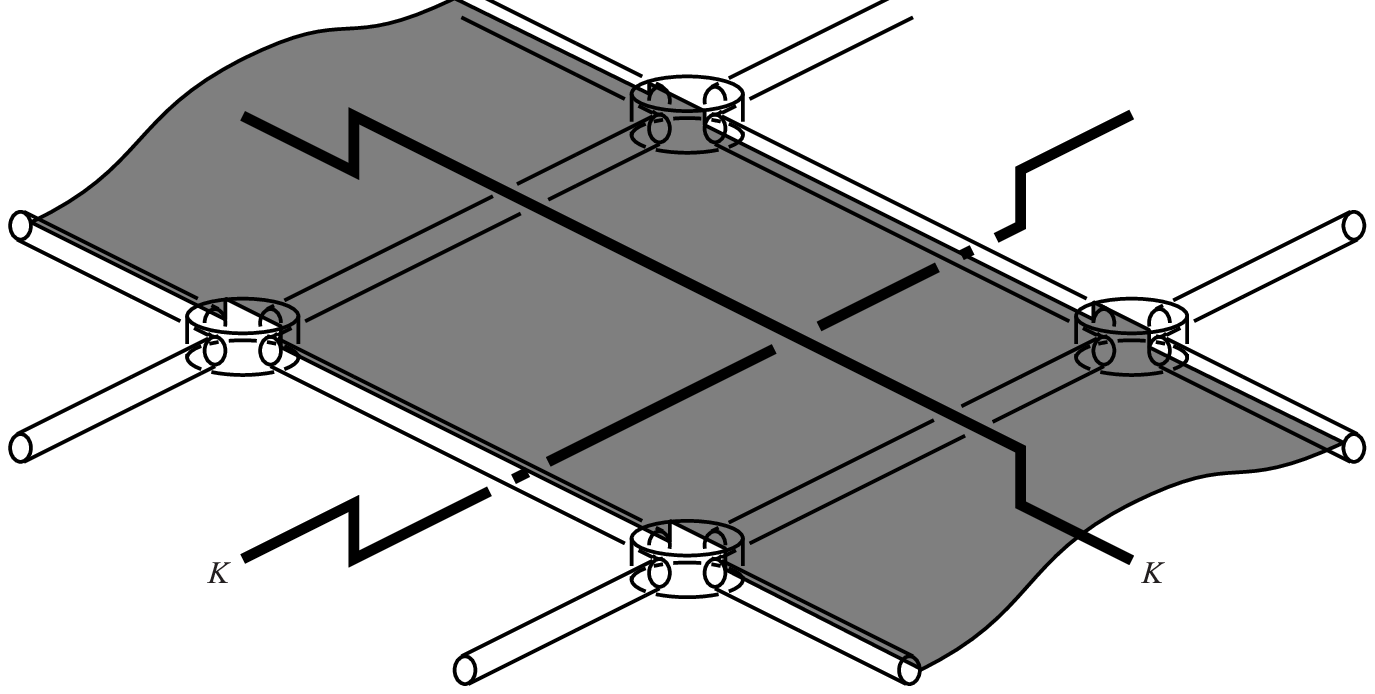}
}
\centerline{Figure 5.}

Finally, there are two 3-handles, one being the 3-ball that lies above
all the handles we have just described, the other being the 3-ball
that lies below these handles.
Thus, we have defined the handle structure ${\cal H}'$. Its underlying space is
clearly the exterior of $K$. 

We need to specify slightly more precisely how the 2-handles run over the
0-handles and 1-handles.
Consider a 0-handle $H_0' = D^2 \times [-1,1]$ of ${\cal H}'$. 
Its intersection with the 1-handles and 2-handles is the regular neighbourhood
of a graph, where the intersection with the 1-handles is a collection
of thickened vertices, and the intersection with the 2-handles forms thickened
edges. The four thickened
vertices are arranged cyclically around the annulus $\partial D^2 \times
[-1,1]$. Thus, we may speak of two thickened
vertices in $\partial H'_0$ as being {\sl opposite} or {\sl adjacent}.
We may arrange that each thickened edge
that runs between adjacent thickened vertices is a thickening of an arc
$\beta \times p$, where $\beta$ is an arc in $\partial D^2$ and $p$ is a point
of $[-1,1]$. We may also arrange that each thickened
edge running between opposite vertices intersects $\partial D^2 \times [-1,1]$
in two thickened vertical arcs, and intersects $D^2 \times \{-1,1 \}$
in a thickened Euclidean geodesic. Thus, it is not hard to see that
the way that the 1-handles and 2-handles of ${\cal H}'$ are attached
to $H'_0$ is as shown in Figure 6. In fact, this
specific arrangement is that of the 0-handle at the bottom of
Figures 2 - 5. 

\vskip 18pt
\centerline{
\epsfxsize=1.8in
\epsfbox{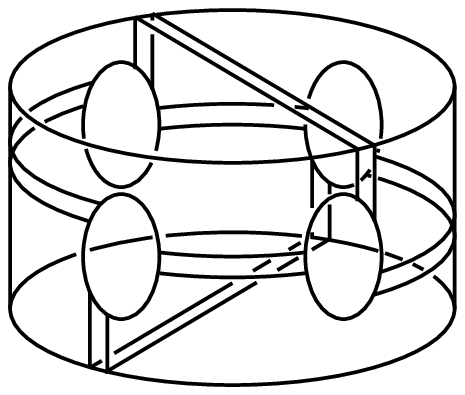}
}
\vskip 6pt
\centerline{Figure 6.}

When the precise embedding of the 0-handle $H'_0$ in $S^3$ is immaterial, we will
usually distort the above picture so that the intersection between $H'_0$ and
the 1-handles and 2-handles is planar, as shown in Figure 7.

\vskip 18pt
\centerline{
\epsfxsize=1.7in
\epsfbox{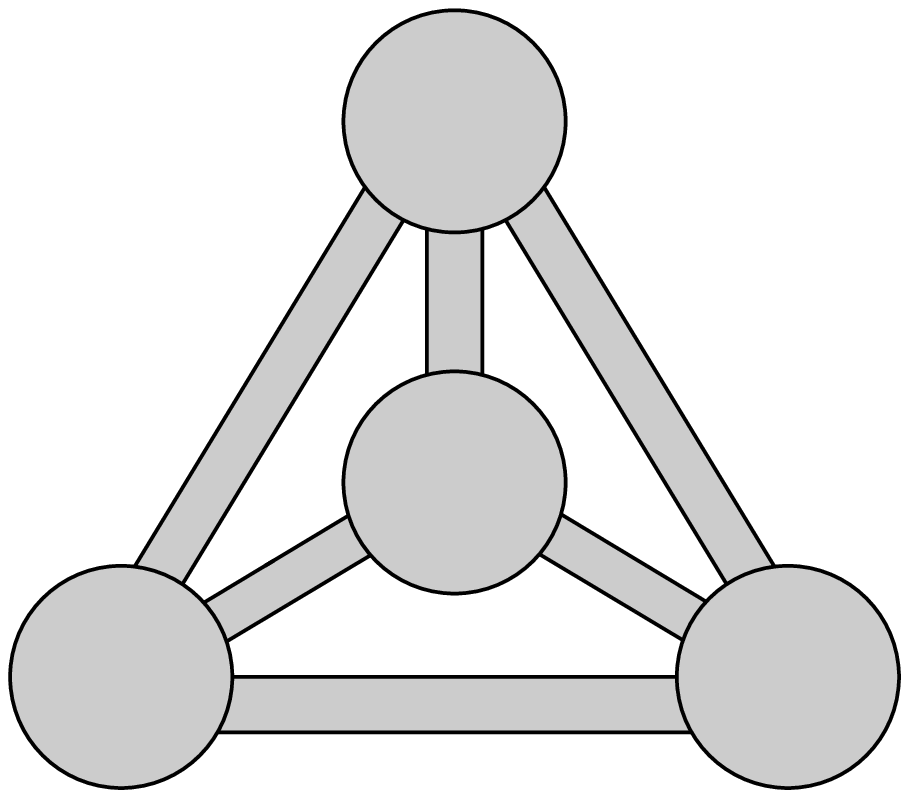}
}
\vskip 6pt
\centerline{Figure 7.}

We now modify ${\cal H}'$ slightly to give the handle
structure ${\cal H}$. Pick a point on $K$ away from the crossings. Running either side of this
point, parallel to $K$, are two 1-handles. Subdivide these, by introducing
a $0$-handle into each. Above and below the knot at this point, there are two
$2$-handles. Subdivide each of these, by introducing a $1$-handle into each, which
runs between the two new 0-handles. (See Figure 8.) Now remove one of these
newly introduced 2-handles above $K$, cancelling it with the 3-handle that
lies above the diagram. Do the same with the 2-handle directly below
it. (See Figure 8.)

\vskip 18pt
\centerline{
\epsfxsize=4.9in
\epsfbox{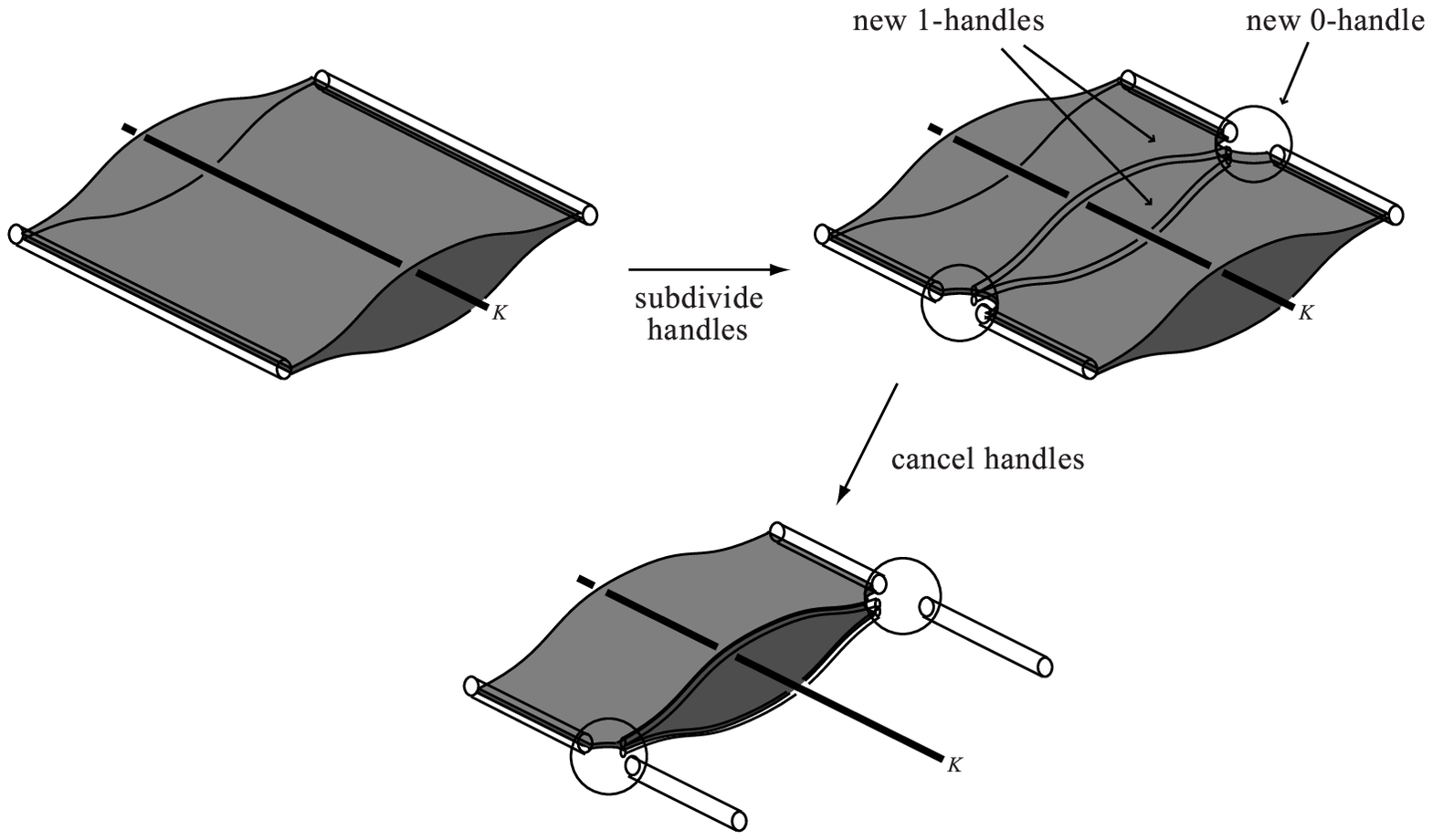}
}
\centerline{Figure 8.}

Let ${\cal H}$ be the resulting handle structure, which we call
a {\sl diagrammatic handle structure}. We term the two new 0-handles of ${\cal H}$
that do not lie in ${\cal H}'$ as {\sl exceptional}. The remaining
0-handles are {\sl unexceptional}. For each $i \in \{ 0,1,2,3 \}$,
let ${\cal H}^i$ be the union of the $i$-handles of ${\cal H}$.
Let ${\cal F}$ be the surface $\partial {\cal H}^0 \cap ({\cal H}^1 \cup {\cal H}^2)$.
As above, this surface may be viewed as the regular neighbourhood of a graph,
with the thickened vertices (denoted ${\cal F}^0$) being ${\cal H}^0 \cap {\cal H}^1$
and the thickened edges (denoted ${\cal F}^1$) being ${\cal H}^0 \cap {\cal H}^2$.
Note that when $H_0$ is an unexceptional 0-handle of ${\cal H}$, then
$H_0 \cap {\cal F}$ is as shown in Figures 6 and 7, but possibly with
some thickened edges removed.

\vskip 18pt
\centerline{\caps 4. Normal surfaces in handle structures}
\vskip 6pt

We now have a diagrammatic handle structure ${\cal H}$ on the exterior
of the knot $K$ arising from the diagram $D$.
When $K$ is a connected sum $K_1 \sharp \dots \sharp K_n$ of non-trivial knots $K_1, \dots, K_n$,
recall that there are associated annuli $A_1, \dots, A_n$
properly embedded in the exterior of $K$. Let $A$ be
their union. A key step in our argument is to place $A$
into normal form with respect to ${\cal H}$. In this section, we recall
what is meant by normal surfaces in a handle structure ${\cal H}$ on
a compact 3-manifold $X$. 
As in Section 3, we denote the union of
the $i$-handles of a handle structure ${\cal H}$ by ${\cal H}^i$. 

\noindent {\bf Convention 4.1.} We will insist throughout this paper that 
any handle structure on a 3-manifold satisfies the following conditions:
\item{(i)} each $i$-handle $D^i \times D^{3-i}$ intersects $\bigcup_{j \leq i-1} {\cal H}^j$
in $\partial D^i \times D^{3-i}$;
\item{(ii)} any two $i$-handles are disjoint;
\item{(iii)} the intersection of any 1-handle
$D^1 \times D^2$ with any 2-handle $D^2 \times D^1$ is of
the form $D^1 \times \alpha$ in $D^1 \times D^2$, where
$\alpha$ is a collection of arcs in $\partial D^2$,
and of the form $\beta \times D^1$ in $D^2 \times D^1$,
where $\beta$ is a collection of arcs in $\partial D^2$;
\item{(iv)} each 2-handle of ${\cal H}$ runs over at least one 1-handle.

\noindent The diagrammatic handle structure constructed in Section 3 satisfies these
requirements.

Let ${\cal F}$ be the surface ${\cal H}^0 \cap ({\cal H}^1 \cup {\cal H}^2)$,
let ${\cal F}^0$ be ${\cal H}^0 \cap {\cal H}^1$, and 
let ${\cal F}^1$ be ${\cal H}^0 \cap {\cal H}^2$.
By the above conditions, ${\cal F}$ is a thickened graph,
where the thickened vertices are ${\cal F}^0$ and the thickened
edges are ${\cal F}^1$.

\noindent {\bf Definition 4.2.} 
We say that a surface $A$ properly embedded in $X$ is {\sl standard}
if
\item{(i)} it intersects each 0-handle in a collection of properly
embedded disjoint discs;
\item{(ii)} it intersects each 1-handle $D^1 \times D^2$ in 
$D^1 \times \beta$, where $\beta$ is a collection of properly
embedded disjoint arcs in $D^2$;
\item{(iii)} it intersects each 2-handle $D^2 \times D^1$ in 
$D^2 \times P$, where $P$ is a collection of points in the
interior of $D^1$;
\item{(iv)} it is disjoint from the 3-handles.

\vskip 18pt
\centerline{
\epsfxsize=4in
\epsfbox{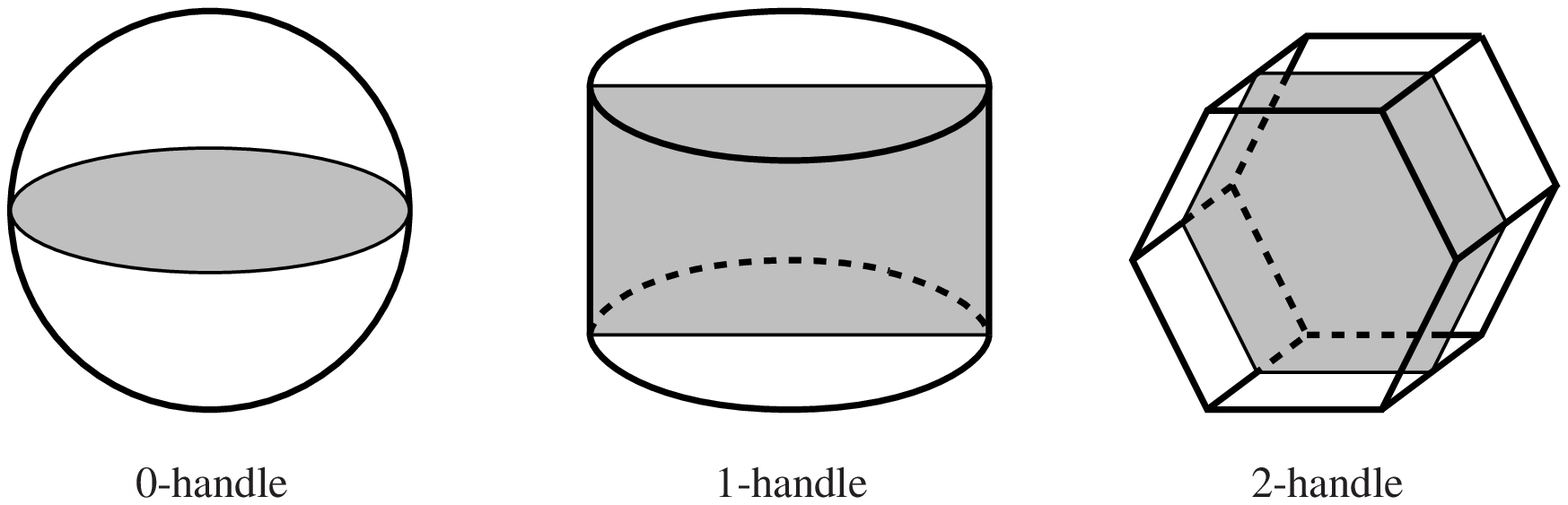}
}
\vskip 6pt
\centerline{Figure 9.}

A standard surface $A$ is termed {\sl normal} if its intersection
with the 0-handles satisfies some conditions, as follows. 

\vfill\eject
\noindent {\bf Definition 4.3.} A disc component $D$ of $A \cap {\cal H}^0$ is 
said to be {\sl normal} if
\item{(i)} $\partial D$ intersects any thickened edge of ${\cal F}$ in
at most one arc;
\item{(ii)} $\partial D$ intersects any component of $\partial {\cal F}^0 - {\cal F}^1$
at most once;
\item{(iii)} $\partial D$ intersects any component of $\partial {\cal H}^0 - {\cal F}$
in at most one arc and no simple closed curves.

\noindent A standard surface that intersects each 0-handle in a disjoint
union of normal discs is said to be {\sl normal}. (See Figure 10.)

\vskip 12pt
\centerline{
\epsfxsize=2.6in
\epsfbox{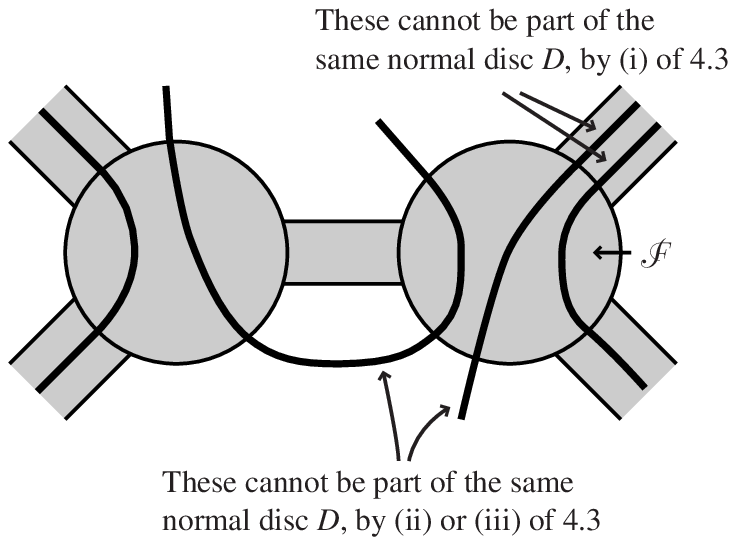}
}
\vskip 6pt
\centerline{Figure 10.}

This is a slightly weaker version of normality than is used by some authors,
for example Definition 3.4.1 in [5]. However, if we had used the definition in [5],
Proposition 4.4 (below) would no longer have held.

If $A$ is a normal surface in $X$, then we also say that a component of
intersection between $A$ and a 1-handle or 2-handle of ${\cal H}$ is a {\sl normal disc}.

Let $H$ be a handle of ${\cal H}$.
Then two normal discs $D$ and $D'$ in $H$ are {\sl normally isotopic} if there an ambient
isotopy, preserving each handle of ${\cal H}$, taking $D$ to $D'$. The discs are then said to be of the
same {\sl normal disc type}. It is a standard fact in normal surface theory that, for each
handle $H$ in a handle structure, there is an upper bound on the
number of normal disc types in $H$, and these disc types are all
constructible. (See p.140 in [5] for example.) 
Indeed, when $H_0$ is an unexceptional 0-handle of the diagrammatic handle structure
and $D$ is a normal disc in $H_0$ that is disjoint from $\partial X$, then $D$ runs over three or
four thickened edges, forming either a {\sl triangle} or {\sl square},
as in Figure 11. However, there is a multitude of different normal disc
types which intersect $\partial X$.

\vskip 18pt
\centerline{
\epsfxsize=3in
\epsfbox{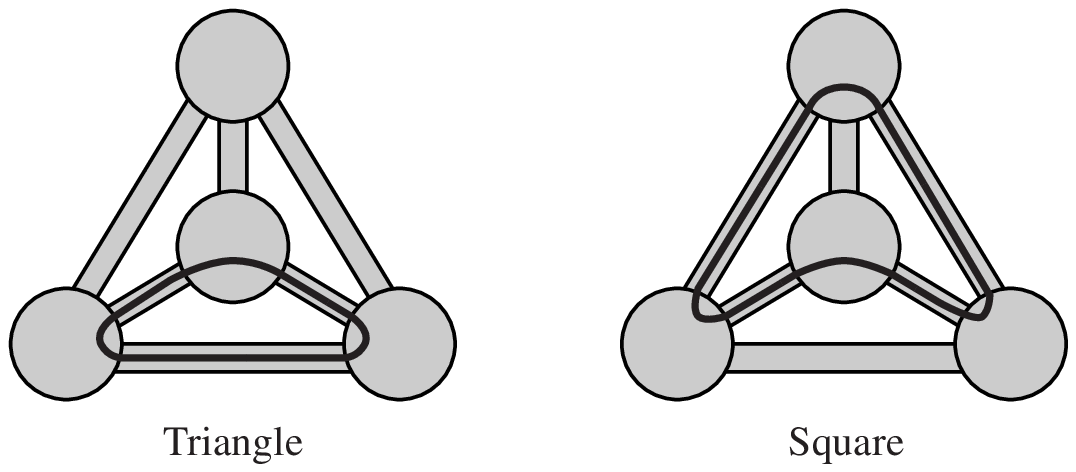}
}
\vskip 6pt
\centerline{Figure 11.}

We now introduce a notion of standard curves in the boundary
of a 3-manifold. Let ${\cal H}$ be a handle structure on a compact
3-manifold $X$. Then a collection of disjoint simple closed curves
in $\partial X$ is {\sl standard} if
\item{(i)} it is disjoint from the 2-handles;
\item{(ii)} it intersects each 1-handle $D^1 \times D^2$ in
$D^1 \times P$, where $P$ is a collection of points in $\partial D^2$;
\item{(iii)} it intersects ${\rm cl}(\partial {\cal H}^0 - {\cal F})$
in a collection of properly embedded arcs.

Note that when $A$ is a normal surface in a handle structure ${\cal H}$
of a 3-manifold $X$, the manifold $M$ obtained by cutting $X$ along
$A$ inherits a handle structure ${\cal H}'$. Note moreover that
the copies of $\partial A$ in $\partial M$ are standard simple closed curves
in ${\cal H}'$.

We now wish to place the annuli $A$ into normal form in the handle structure
${\cal H}$ on $X$, the exterior of $K$. We first ambient isotope $\partial A$
so that it runs over the two exceptional 0-handles and the two 1-handles
that run between them, as shown in Figure 12. Note that $\partial A$
is then standard in $\partial X$.
Recall that $A$ divides $X$ into 3-manifolds $X_1, \dots, X_n$
and $Y$, where each $X_i$ is a copy of the exterior of $K_i$.
We may also ensure that for each $i \in \{ 1, \dots, n \}$,
$\partial X_i \cap \partial X$ also lies in the exceptional 0-handles
and the 1-handles that run between them.

\vskip 12pt
\centerline{
\epsfxsize=2in
\epsfbox{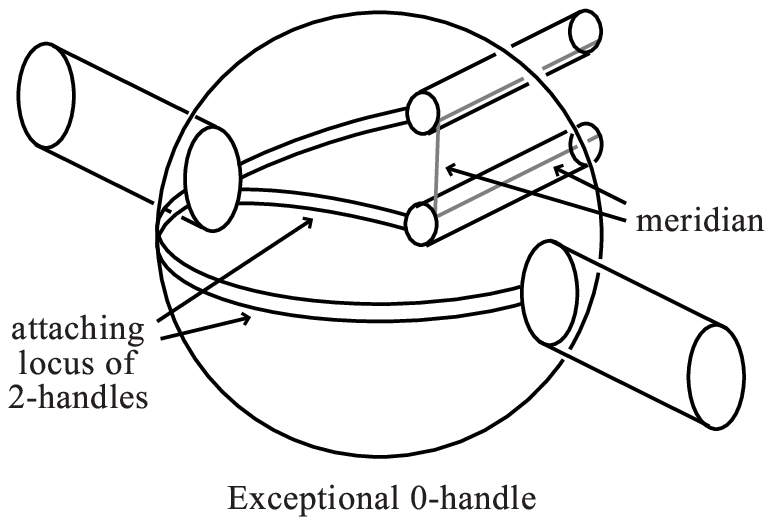}
}
\vskip 6pt
\centerline{Figure 12.}
\vfill\eject

We now apply the following proposition, which is a slight variant of
a well-known result in normal surface theory.

\noindent {\bf Proposition 4.4.} {\sl Let ${\cal H}$ be a handle structure on
a compact irreducible 3-manifold $X$. Let $A$ be a properly embedded, incompressible,
boundary-incompressible surface in $X$, with no 2-sphere components. Suppose that
each component of $\partial A$ is standard in $\partial X$ and intersects
each component of $\partial X \cap {\cal H}^0$ and $\partial X \cap {\cal H}^1$
in at most one arc and no simple closed curves. Then
there is an ambient isotopy, supported in the interior of $X$, taking $A$
into normal form.}

\noindent {\sl Proof.} This is fairly routine, and follows the proof of
Theorem 3.4.7 in [5] for example. So, we will
only sketch the argument. We may first ambient isotope $A$ so that it 
misses the 3-handles of ${\cal H}$ and respects the
product structure on the 1-handles and 2-handles. Suppose that some component
of intersection between $A$ and a 0-handle $H_0$ is not a disc.
Then $A \cap H_0$ is compressible, via a compression disc $D_1$.
Since $A$ is incompressible, $\partial D_1$ bounds a disc $D_2$ in $A$.
As $X$ is irreducible, $D_1 \cup D_2$ bounds a 3-ball, and we may ambient
isotope $D_2$ onto $D_1$. This either reduces $|A \cap {\cal H}^2|$, or
it leaves it unchanged and reduces $|A \cap {\cal H}^1|$. So, we may
assume that $A$ intersects each 0-handle in a collection of discs. If
a component of intersection between $A$ and a 1-handle is not a disc, then
it is an annulus, which forms part of a 2-sphere component of $A$, contrary
to assumption. Thus, $A$ is now standard. Consider a component $D$ of $A \cap {\cal H}^0$.
If this intersects some thickened edge of ${\cal F}$ more than once, then
there is an ambient isotopy, which reduces $|A \cap {\cal H}^2|$. So,
we may assume that (i) in Definition 4.3 (the definition of normality) holds. Suppose
that $\partial D$ intersects a component of $\partial {\cal H}^0 - {\cal F}$
in more than one arc. These two arcs may be joined by an arc $\alpha$ in
$\partial {\cal H}^0 - {\cal F}$. The endpoints of $\alpha$ may be joined
by a properly embedded arc $\beta$ in $D$. By choosing $D$ suitably,
we may ensure that the interior of $\alpha$ is disjoint from $A$.
Then $\alpha \cup \beta$ bounds a disc $D'$ in ${\cal H}_0$ such that
$D' \cap \partial {\cal H}^0 = \alpha$ and $D' \cap A = \beta$.
Now, $A$ is boundary-incompressible, and so $\beta$ separates $A$ into
two components, one of which is a disc. In particular,
the endpoints of $\alpha$ lie in the
same component of $\partial A$. Hence, this component of $\partial A$ intersects a component of
$\partial X \cap {\cal H}^0$ in more than one arc, which is contrary to hypothesis.
Also, $\partial D$ cannot intersect $\partial {\cal H}^0 - {\cal F}$
in a simple closed curve, by hypothesis.
Thus, (iii) in Definition 4.3 is verified. Finally, $\partial D$ cannot intersect $\partial {\cal F}^0 - {\cal F}^1$
more than once, since this would imply that a component of $\partial A$
runs over a component of $\partial X \cap {\cal H}^1$ more than once. Thus, $D$
satisfies (ii) of Definition 4.3, and so $A$ is normal. $\square$

\vfill\eject
\centerline{\caps 5. Generalised parallelity bundles}
\vskip 6pt

Recall that we are going to cut the exterior of $K$
along the annuli $A$. The result will be the disjoint union of
3-manifolds $X_1, \dots, X_n$ and $Y$, where each $X_i$ is homeomorphic to
the exterior of $K_i$. We will choose on the boundary of each $X_i$
(after $X_1 \cup \dots \cup X_n$ has been ambient isotoped)
a simple closed curve $C_i$ which hits a meridian of $K_i$ just once.
Then, the diagram $D'$ for $K_1 \sqcup \dots \sqcup K_n$ will be the projection of $C_1 \cup \dots \cup C_n$.
Our goal is to restrict the number of crossings of $D'$.

Now, the boundary of $X_1 \cup \dots \cup X_n$ is partitioned into two subsurfaces: a
copy of $A$ and parts of $\partial N(K)$. Let $S$ be the former
surface. The parts of $\partial N(K)$ in $X_1 \cup \dots \cup X_n$ form a
reasonably controlled subsurface of $S^3$. However, the annuli $S$
may be complicated. They consist of normal discs (one for each
normal disc of $A$), but $A$ may be made up of many normal discs.
One might hope to prove
Theorem 1.1 by bounding the total number
of normal discs of $A$. However, to prove Theorem 1.1
in this way, this bound would need to be
{\sl linear} in the number of crossings
of $D$. Now, there are known bounds on the
number of normal discs of certain surfaces in
handle structures, for example, Lemma 3.2 of [2].
But these are exponential in the number of
0-handles of ${\cal H}$. It seems unlikely that
one can achieve a linear bound in general.
Thus, a new approach is required. Suppose that
$A$ intersects a handle of ${\cal H}$ in many normal discs. Then
many of these must be normally parallel. The region
between two adjacent parallel normal discs of $A$ is a product
$D \times I$, where $(D \times I) \cap A = D \times \partial I$.
These parallelity regions combine to form $I$-bundles
embedded in the exterior of $A$. In this section,
we consider a generalisation of this structure, known as a `generalised
parallelity bundle'.

Let $M$ be a compact orientable 3-manifold with a handle structure ${\cal H}$,
and let $S$ be a subsurface of $\partial M$ such that $\partial S$ is
standardly embedded in $\partial M$.
We then say that ${\cal H}$ is a {\sl handle
structure} for the pair $(M,S)$.
The main example we will consider is where $M = X_1 \cup \dots \cup X_n$,
and $S$ is the copy of $A$ in $M$.

\noindent {\bf Definition 5.1.}
A handle $H$ of ${\cal H}$ is a {\sl parallelity handle} if
it admits a product structure $D^2 \times I$ such that
\item{(i)} $D^2 \times \partial I = H \cap S$;
\item{(ii)} each component of ${\cal F}^0 \cap H$ and ${\cal F}^1 \cap H$
is $\beta \times I$, for a subset $\beta$ of $\partial D^2$.

We will typically view the product structure $D^2 \times I$
as an $I$-bundle over $D^2$.

The main example of a parallelity handle arises when $M$ is obtained by cutting a 3-manifold $X$
along a normal surface $A$, and where $S$ is the copies of $A$ in $M$. Then, if $A$ contains 
two normal discs in a handle that are normally 
parallel and adjacent, the space between them becomes a parallelity
handle in $M$. See Figure 13.

\vskip 18pt
\centerline{
\epsfxsize=3.5in
\epsfbox{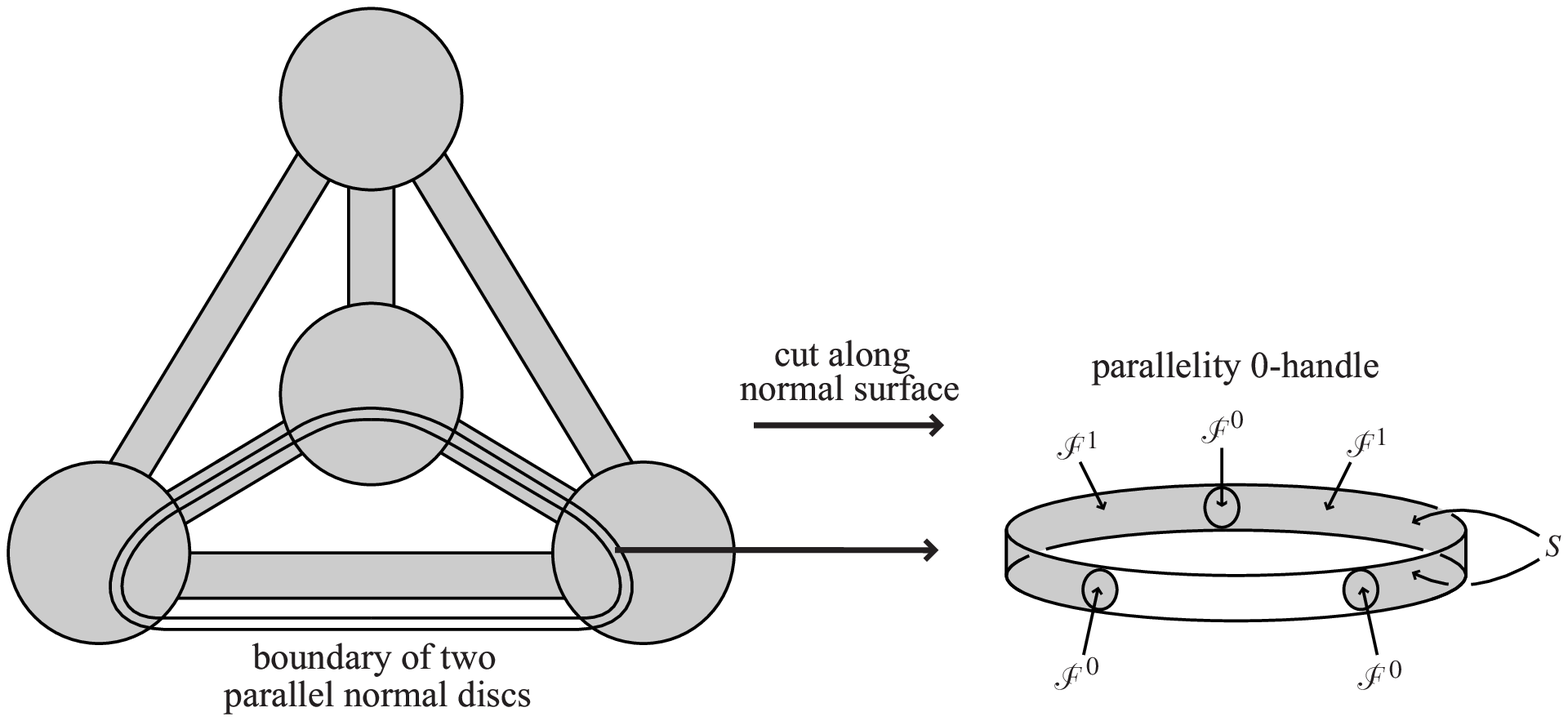}
}
\vskip 6pt
\centerline{Figure 13.}

We now collate some facts about parallelity handles.

A 3-handle can never be a parallelity handle because 3-handles are disjoint
from $\partial M$, and hence disjoint from $S$.

A 2-handle $D^2 \times D^1$ is a parallelity handle if and only if
$(D^2 \times D^1) \cap S = D^2 \times \partial D^1$. In this case,
the two product structures that the handle has, one from the
fact that it is a parallelity handle, the other from the fact
that is a 2-handle, can be made to coincide.

A 1-handle $D^1 \times D^2$ is a parallelity handle if and only if
$(D^1 \times D^2) \cap S$ is two discs and each component of
$(D^1 \times \partial D^2) - S$ lies entirely in $\partial M$ or entirely in
${\cal H}^2$. In this case, the $I$-bundle
structure on the 1-handle can be made to respect its
structure as a product $D^1 \times D^2$, in the sense that $D^2$
inherits a structure as $I \times I$, so that fibres in the
$I$-bundle are of the form $p_1 \times p_2 \times I$, for $p_1 \in D^1$
and $p_2 \in I$. 

When two parallelity handles are incident, we will see that their $I$-bundle structures
can be made to coincide along their intersection. So, the union of the parallelity
handles forms an $I$-bundle over a surface $F$, say. (See Lemma 5.3.)
It will be technically convenient to consider enlargements of
such structures. These will still be an $I$-bundle over a surface
$F$, and near the $I$-bundle over $\partial F$, they will be a union
of parallelity handles, but elsewhere need not be. The precise definition
is as follows.

\noindent {\bf Definition 5.2.}
Let ${\cal H}$ be a handle structure for the pair $(M,S)$.
A {\sl generalised parallelity bundle} ${\cal B}$ is a 3-dimensional
submanifold of $M$ such that
\item{(i)} ${\cal B}$ is an $I$-bundle over 
a compact surface $F$;
\item{(ii)} the $\partial I$-bundle is ${\cal B} \cap S$;
\item{(iii)} ${\cal B}$ is a union of handles of ${\cal H}$;
\item{(iv)} any handle in ${\cal B}$ that intersects the $I$-bundle
over $\partial F$ is a parallelity handle, where $I$-bundle
structure on the parallelity handle agrees with the $I$-bundle structure of
${\cal B}$;
\item{(v)} ${\rm cl}(M - {\cal B})$ inherits
a handle structure.

\noindent The $I$-bundle over $\partial F$ is termed
the {\sl vertical boundary} of ${\cal B}$,
and the $\partial I$-bundle over $F$ is called
the {\sl horizontal boundary}.

Note that a single 2-handle $D^2 \times D^1$ such
that $(D^2 \times D^1) \cap S = D^2 \times \partial D^1$ is a
generalised parallelity bundle. An example of a slightly more
complicated generalised parallelity bundle is shown in Figure 14.
It is composed of two parallelity 2-handles and a parallelity
1-handle.

\vskip 18pt
\centerline{
\epsfxsize=2.6in
\epsfbox{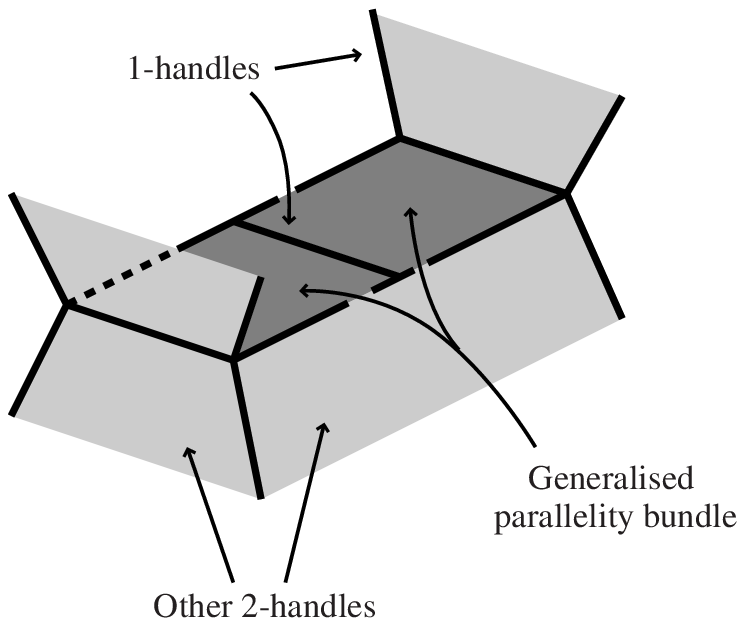}
}
\vskip 6pt
\centerline{Figure 14.}

The point of the
definition is that generalised parallelity bundles behave in many
ways like 2-handles, in that the remaining handles form a handle
structure, onto which the generalised parallelity bundle is
attached. However, note that the vertical boundary of a generalised parallelity bundle 
need not be properly embedded in $M$. This is because the vertical
boundary of a parallelity handle may intersect $\partial M$ in its
interior. But, the intersection between the vertical boundary
and $\partial M$ is a union of fibres in the $I$-bundle.

The following lemma gives an important example of a generalised
parallelity bundle.

\noindent {\bf Lemma 5.3.} {\sl The union of the parallelity handles
is a generalised parallelity bundle.}

\noindent {\sl Proof.} By definition, each parallelity handle has the
structure of an $I$-bundle. We claim that these structures can be chosen
so that they coincide on the intersection of any two parallelity handles. Hence, the
union ${\cal B}$ of the parallelity handles will inherit an $I$-bundle structure.

Consider first the intersection of a parallelity 0-handle $H_0$ and
a parallelity 1-handle $H_1$. These intersect along components
of ${\cal F}^0$. By condition (ii) in the Definition 5.1 (the definition of a parallelity
handle), each such component of ${\cal F}^0$ inherits a product structure
$\beta_0 \times I$ from $H_0$ and a product structure $\beta_1 \times I$
from $H_1$. By condition (i) in Definition 5.1 applied twice,
$$\beta_0 \times \partial I = (\beta_0 \times I) \cap S = (\beta_1 \times I) \cap S
= \beta_1 \times \partial I.$$
Hence, these product structures can be made to coincide. A similar
argument applies to the intersection of a parallelity 0-handle and
a parallelity 2-handle, but with the role of ${\cal F}^0$ replaced by
${\cal F}^1$. Finally consider the intersection of a parallelity
1-handle $H_1$ and a parallelity 2-handle $H_2$. 
Now, $H_1 \cap {\cal F}^1$ is a disjoint union of fibres, by condition (ii) in
the Definition 5.1. Thus, the $I$-bundle
structures of $H_1$ and $H_2$ agree along $H_1 \cap H_2 \cap {\cal H}^0$.
Since the $I$-bundle structures respect the product structures
on $H_1$ and $H_2$, we see that they agree along all of $H_1 \cap H_2$.

Thus, conditions (i) - (iv) in Definition 5.2 (the definition of a generalised parallelity
bundle) follow immediately. We must check condition (v) in Definition 5.2, which asserts that
${\rm cl}(M - {\cal B})$ inherits a handle structure. The only way that this
might fail is if a $j$-handle of ${\rm cl}(M - {\cal B})$ is incident
to an $i$-handle of ${\cal B}$, for $j > i$. Thus, we must check that
if an $i$-handle is a parallelity handle, then
so is any $j$-handle to which it is incident, for $j > i$. Let us consider
when the $i$-handle is a 0-handle $H_0$. Then, by definition of
the parallelity structure on $H_0$, each component of $H_0 \cap {\cal F}^0$
and $H_0 \cap {\cal F}^1$ inherits an $I$-bundle structure, which
therefore extends over any incident 1-handle or 2-handle, making it
a parallelity handle. Note also that a parallelity 0-handle
is not incident to any 3-handles. So, the claim holds for $i = 0$. 
Let us now consider the case where $i=1$, and let $H_1 = D^1 \times D^2$ be a parallelity
1-handle. Then $H_1 \cap S$ is two discs, and each component of $(D^1 \times \partial D^2) - S$
lies entirely in $\partial M$ or entirely in ${\cal H}^2$. Thus, $H_1$
is disjoint from the 3-handles. Also, any 2-handle $D^2 \times D^1$ to which it
is incident has both components of $D^2 \times \partial D^1$ lying in $S$.
So, it is a parallelity handle, as required. Finally, the case where
$i=2$ follows from the observation that a parallelity 2-handle is
disjoint from the 3-handles. This proves the claim.
It is now clear that conditions (i) - (iv) in Convention 4.1 hold for
${\rm cl}(M - {\cal B})$, and hence it inherits a handle structure.
$\square$

Suppose that $M$ is irreducible and $S$ is incompressible.
Our aim now is to construct a handle structure
on $(M,S)$ containing a generalised parallelity bundle that satisfies the following
two conditions:
\item{(i)} it contains every parallelity handle;
\item{(ii)} its horizontal boundary is incompressible in $M$.

\noindent This will be achieved via the following procedure for simplifying a
handle structure ${\cal H}$ of $(M,S)$.

\noindent {\bf Definition 5.4.}
Let $G$ be an annulus properly embedded in $M$, with boundary in $S$.
Suppose that there is an annulus $G'$ in $\partial M$ such that $\partial G = \partial G'$.
Suppose also that $G \cup G'$ bounds a 3-manifold $P$ such that
\item{(i)} either $P$ is a parallelity region between $G$ and $G'$, or $P$ lies
in a 3-ball;
\item{(ii)} $P$ is a non-empty union of handles;
\item{(iii)} ${\rm cl}(M - P)$ inherits
a handle structure from ${\cal H}$;
\item{(iv)} any parallelity handle of ${\cal H}$ that intersects
$P$ lies in $P$;
\item{(v)} $G$ is a vertical boundary component of a generalised
parallelity bundle lying in $P$;
\item{(vi)} $G' \cap (\partial M - S)$ is either empty or a
regular neighbourhood of a core curve of the annulus $G'$.

\noindent Removing the interiors of $P$
and $G'$ from $M$ is called an {\sl annular
simplification}. Note that the resulting 3-manifold
$M'$ is homeomorphic to $M$, even though $P$ may
be homeomorphic to the exterior of a non-trivial knot
when it lies in a 3-ball. (See Figure 15.)
The boundary of $M'$ inherits a copy of $S$, which we denote
by $S'$, as follows. We set $S' \cap \partial M$ to be
$S \cap \partial M'$. When $G' \cap (\partial M - S)$ is empty,
we declare that $\partial M' - \partial M$ lies in $S'$.
When $G' \cap (\partial M - S)$ is a single annulus,
we declare that $\partial M' - \partial M$ is disjoint from $S'$.
Thus, $(M',S')$ is homeomorphic to $(M,S)$. Moreover, when $M$ is embedded 
within a bigger closed 3-manifold, then
$(M',S')$ is ambient isotopic to $(M,S)$.

\vskip 18pt
\centerline{
\epsfxsize=3.0in
\epsfbox{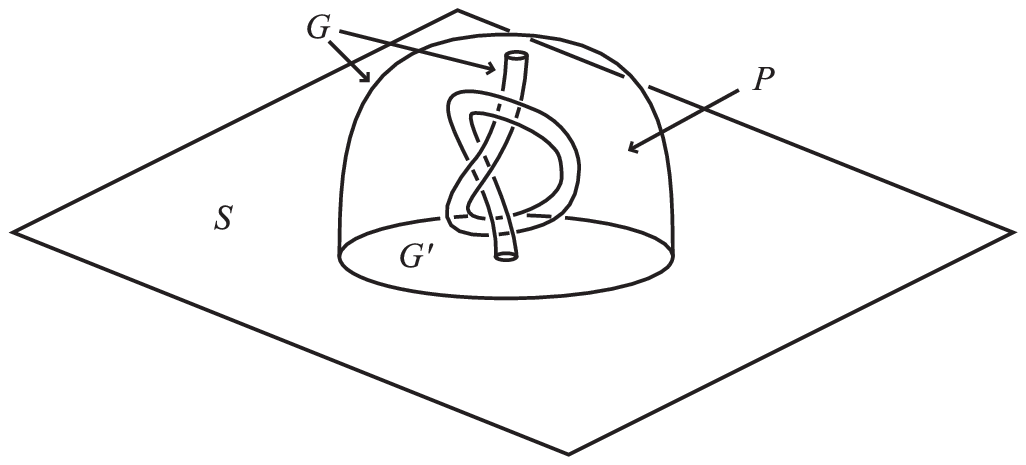}
}
\vskip 6pt
\centerline{Figure 15.}

\noindent {\bf Lemma 5.5.} {\sl Let ${\cal H}$ be a handle structure for the
pair $(M,S)$. Let ${\cal H}'$ be a handle structure obtained from ${\cal H}$ by
an annular simplification. Then any parallelity handle for ${\cal H}$ that
lies in ${\cal H}'$ is a parallelity handle for ${\cal H}'$.}

\noindent {\sl Proof.} Let $H$ be a parallelity handle for ${\cal H}$.
Let $P$ be the 3-manifold, the interior of which is removed in the annular
simplification. If $H$ intersects $P$, then it lies in $P$ by 
condition (iv) in Definition 5.4 (the definition of an annular simplification), and hence it does not
lie in ${\cal H}'$. So, if $H$ lies in ${\cal H}'$, then it
does not intersect $P$, and so it was not modified in the annular
simplification. Thus, it is a parallelity handle for ${\cal H}'$.
$\square$

A generalised parallelity bundle is {\sl maximal} if it is not strictly contained within another
generalised parallelity bundle. (Thus, it is maximal with respect to the partial
order of inclusion, where the inclusion does not necessarily respect the
bundle structures.) Note that if a generalised parallelity bundle ${\cal B}$
is maximal, then any parallelity handle that intersects ${\cal B}$
must lie in ${\cal B}$.

\noindent {\bf Proposition 5.6.} {\sl Let $M$ be a compact
orientable irreducible 3-manifold with a handle structure ${\cal H}$.
Let $S$ be an incompressible subsurface of $\partial M$, such that
$\partial S$ is standard in $\partial M$.
Suppose that ${\cal H}$ admits no annular simplification.
Let ${\cal B}$ be any maximal generalised parallelity bundle in ${\cal H}$.
Then the horizontal boundary of ${\cal B}$ is incompressible.}

\noindent {\sl Proof.} Let ${\cal B}'$ be those
components of ${\cal B}$ that are not $I$-bundles
over discs. It clearly suffices to show
that the horizontal boundary of ${\cal B}'$
is incompressible.

We claim that it suffices to show that the vertical
boundary of ${\cal B}'$ is incompressible. For, if the vertical boundary
were incompressible, then any compression disc for the
horizontal boundary could be isotoped off the vertical
boundary. Hence, it would lie entirely in the generalised parallelity bundle.
But the horizontal boundary
of an $I$-bundle is incompressible in the $I$-bundle.
Thus, the horizontal
boundary of ${\cal B}'$ is incompressible
if the vertical boundary is.

Consider therefore a compression disc $D$ for the 
vertical boundary of ${\cal B}'$. Let $V$ be
the vertical boundary component containing $\partial D$.
By the definition of ${\cal B}'$, $D$ does
not lie entirely in ${\cal B}'$. Its interior
is disjoint from ${\cal B}'$ (by the definition
of a compression disc), but it may intersect
${\cal B} - {\cal B}'$. Note that $V$ is properly embedded in $M$,
since the interior of $D$ lies on one side of it,
and a component of ${\cal B}'$ lies on the other side.

Now, $V$ compresses along $D$ to give two discs $D'_1$ and $D'_2$
embedded in $M$, with boundary in $S$.
Since $S$ is incompressible and $M$ is irreducible,
$D'_1$ and $D'_2$ are parallel to discs
$D_1$ and $D_2$ in $S$, via 3-balls $P_1$
and $P_2$. There are two cases to consider:
where $P_1$ and $P_2$ are disjoint and where they
are nested.

Let us suppose first that they are disjoint. Then,
$V \cup D_1 \cup D_2$ bounds a 3-ball $B$. Since the interior of $D$ is
disjoint from ${\cal B}'$, this ball $B$ does not lie in ${\cal B}'$. So,
we may extend the $I$-bundle structure of ${\cal B} - B$
over $B$, contradicting the maximality of 
${\cal B}$. See Figure 16. Note that, here, we are using the fact
that generalised parallelity bundles need not consist
solely of parallelity handles.

\vskip 18pt
\centerline{
\epsfxsize=4in
\epsfbox{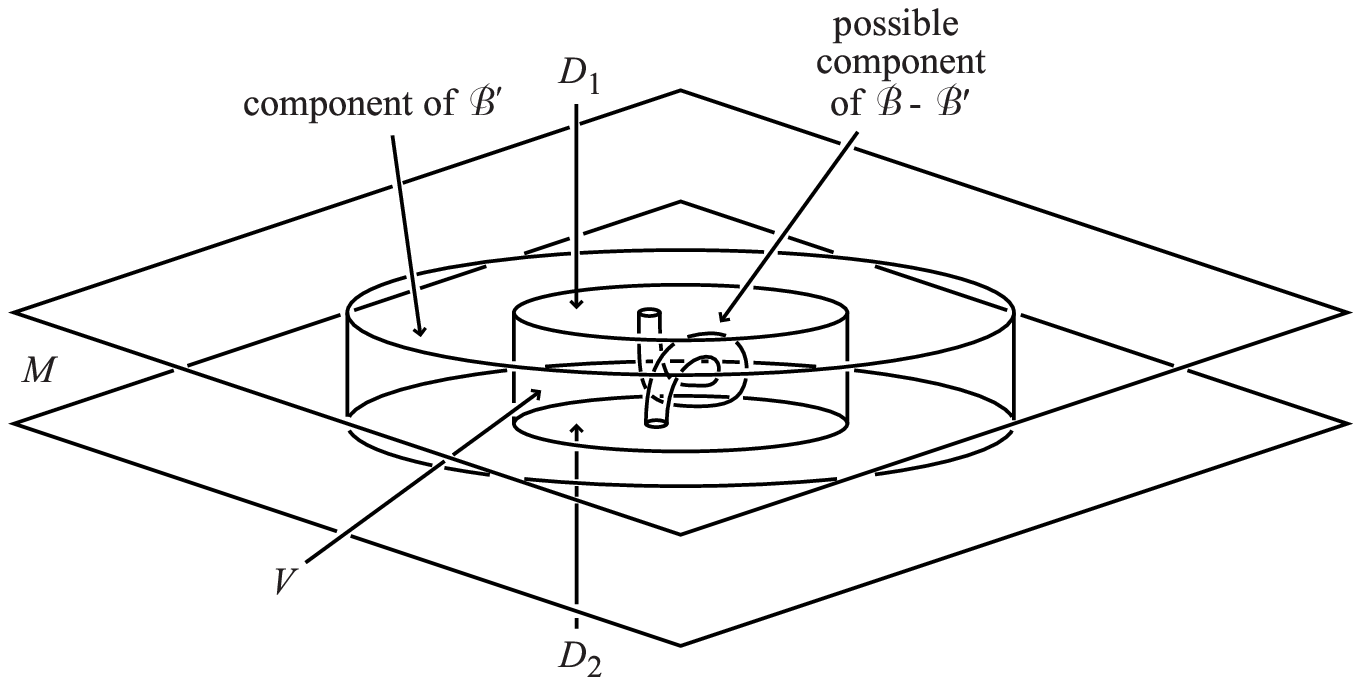}
}
\vskip 6pt
\centerline{Figure 16.}

Let us now suppose that $P_1$ and $P_2$ are nested;
say that $P_2$ lies in $P_1$.  See Figure 17. Let $G'$ be $D_1 - {\rm int}(D_2)$.
Then, $G'$ is an annulus in $S$ such that $\partial V = \partial G'$.
Let $P$ be the 3-manifold bounded by $V \cup G'$.
This lies in the 3-ball $P_1$.
By (v) in Definition 5.2 (the definition of a generalised parallelity bundle), ${\rm cl}(M - P)$
inherits a handle structure. Note also that, by the maximality
of ${\cal B}$, any parallelity handle that intersects
$P$ lies in $P$. So, ${\cal H}$ admits an annular
simplification, which is a contradiction. $\square$

\vskip 18pt
\centerline{
\epsfxsize=3.5in
\epsfbox{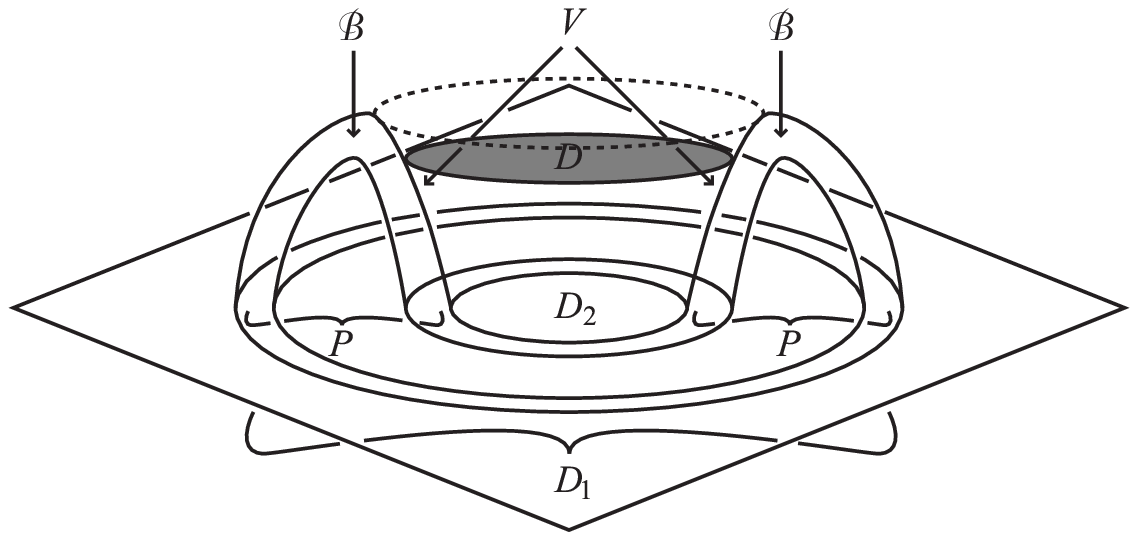}
}
\vskip 6pt
\centerline{Figure 17.}

\noindent {\bf Corollary 5.7.} {\sl Let $M$, $S$ and ${\cal H}$ be as in
Proposition 5.6. Then there is a generalised parallelity bundle
${\cal B}$ such that
\item{(i)} ${\cal B}$ contains every parallelity handle;
\item{(ii)} the horizontal boundary of ${\cal B}$ is incompressible.}

\noindent {\sl Proof.} By Lemma 5.3, the union of the parallelity handles
is a generalised parallelity bundle. Enlarge this to a maximal
generalised parallelity bundle ${\cal B}$. By Proposition 5.6,
its horizontal boundary is incompressible. $\square$

\noindent {\bf Proposition 5.8.} {\sl Let $M$ be a compact, orientable,
irreducible 3-manifold, with boundary a collection of incompressible tori.
Let $S$ be a subsurface of $\partial M$ such that the intersection
of $S$ with each component of $\partial M$ is either empty
or a single incompressible annulus. Suppose that any incompressible
annulus properly embedded in $M$ with boundary in $S$ is boundary parallel.
Let ${\cal H}$ be a handle
structure for $(M,S)$ that admits no annular simplifications.
Then there is a generalised parallelity bundle
${\cal B}$ such that
\item{(i)} ${\cal B}$ contains every parallelity handle;
\item{(ii)} ${\cal B}$ is a collection of $I$-bundles over discs.}

\noindent {\sl Proof.} Let ${\cal B}$ be the generalised parallelity
bundle provided by Corollary 5.7.
Let ${\cal B}'$ be the union of the components of ${\cal B}$ that are not $I$-bundles over
discs. Its horizontal boundary is a subsurface of
$S$, which is a collection of annuli. The only connected compact incompressible subsurface
of an annulus is an annulus or disc. Therefore, each component of ${\cal B}'$ is an $I$-bundle
over an annulus or Mobius band. The vertical boundary
components of ${\cal B}'$ are incompressible annuli, with
boundary curves in $S$. Thus, by assumption, each
such vertical boundary component is boundary-parallel in $M$.

We claim that no component of ${\cal B}'$ is an $I$-bundle over a Mobius band. 
Let $V$ be the vertical boundary of such a component $B$ of ${\cal B}'$. 
Now, $V$ is boundary parallel in $M$, via a parallelity region $P$. The
interior of $P$ is disjoint from $B$, since $B$ is an
$I$-bundle over a Mobius band. Hence, $M$ is homeomorphic
to $B$, which is a solid torus. But, this is a contradiction,
because we have assumed that $\partial M$ is incompressible. This proves the claim.

Consider a component $V$ of the vertical boundary of ${\cal B}'$, 
and let $P$ be the parallelity region between $V$ and an annulus in
$\partial M$. The component of ${\cal B}'$
containing $V$ is an $I$-bundle over an annulus,
and so its vertical boundary components are
parallel. So, by changing the choice of $V$ if necessary,
we may assume that $P$ contains this
component of ${\cal B}'$. Then $V$ must be properly embedded in $M$,
for otherwise one could find a compression disc for
$\partial M$ in $P$.
Hence, if one were to remove the interiors of $P$ and $\partial P - V$
from ${\cal H}$, this would be an annular simplification, contrary
to hypothesis. So, ${\cal B}'$ must be empty, and
therefore ${\cal B}$ is a collection of $I$-bundles over
discs, as required. $\square$

\vskip 18pt
\centerline{\caps 6. Proof of the main theorem}
\vskip 6pt

In this section, we will complete the proof of Theorem 1.1. 
Suppose that $K$ is a connected sum of oriented knots $K_1, \dots, K_n$.

\vfill\eject
\noindent {\caps 6.1. We may assume that each $K_i$ is prime and non-trivial}
\vskip 6pt

Express each $K_i$ as a connected sum of prime knots
$K_{i,1}, \dots, K_{i,m(i)}$. Suppose that we could
prove the theorem in the case where each summand is
prime. Then we would have the inequality
$$c(K) \geq {\sum_{i=1}^n \sum_{j = 1}^{m(i)} c(K_{i,j}) \over 152}.$$
But the trivial inequality for connected sums gives that
$$\sum_{j=1}^{m(i)} c(K_{i,j}) \geq c(K_i),$$
and so this would imply that
$$c(K) \geq {\sum_{i=1}^n c(K_i) \over 152}.$$
Thus, it suffices to consider the case where each $K_i$ is
prime. We may also clearly assume that each $K_i$ is non-trivial.

\vskip 6pt
\noindent {\caps 6.2. Handle structures and normal surfaces}
\vskip 6pt

Let $D$ be a diagram of $K$ with minimal crossing number.
Our aim is to construct a diagram $D'$ for the distant
union $K_1 \sqcup \dots \sqcup K_n$, such that $c(D') \leq 152 \ c(D)$.
This will prove the theorem.

Let $X$ be the exterior of $K$. Give $X$ the diagrammatic handle structure 
described in Section 3. Recall that the expression of $K$ as a 
connected sum $K_1 \sharp \dots \sharp K_n$
specifies a collection of annuli $A_1, \dots, A_n$ properly
embedded in $X$, as shown in Figure 1. Let $A$ be their union.
We first perform the isotopy of $\partial A$ that is described in Section 4
just before Proposition 4.4.
It then intersects only the exceptional 0-handles and the 1-handles that
run between them. (See Figure 12.) Then, using Proposition 4.4, we ambient isotope $A$, keeping
$\partial A$ fixed, taking it to a normal surface.

Let $X_1 \cup \dots \cup X_n \cup Y$ be the result of cutting $X$ along $A$.
Let $M$ be $X_1 \cup \dots \cup X_n$ and let $S$ be the
copy of $A$ in $M$. Let ${\cal H}$ be the handle structure
that $M$ inherits. Note that $\partial S$ is standard in $\partial M$.
Thus, ${\cal H}$ is a handle structure for the pair $(M,S)$.

\vfill\eject
\noindent {\caps 6.3. Applying annular simplifications}
\vskip 6pt

We claim that all of $\partial M \cap \partial X$ lies in the
parallelity handles of ${\cal H}$. Note that there is only one possibility, up to normal isotopy, for
the normal discs of $A$ in the exceptional 0-handles. (See Figure 18.) They are
therefore normally parallel. Moreover, in the isotopy of $\partial A$
in Section 4, we arranged that $\partial M \cap \partial X$ lies
in the exceptional 0-handles and the 1-handles that run between them.
The claim now follows.

\vskip 12pt
\centerline{
\epsfxsize=2in
\epsfbox{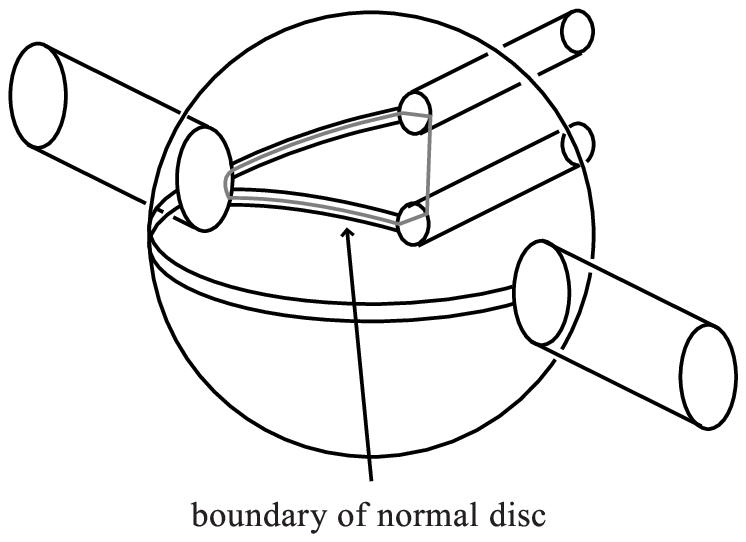}
}
\vskip 6pt
\centerline{Figure 18.}

Let $R$ be the union of the parallelity handles in $M$,
and let $R'$ be the union of the components of $R$ that are
incident to $\partial X$.

We apply as many annular simplifications to ${\cal H}$ as possible,
giving a handle structure ${\cal H}'$ on a pair $(M',S')$ ambient
isotopic to $(M,S)$. Denote the components of $M'$ by
$X'_1, \dots, X'_n$, and let $S'_i = S' \cap \partial X'_i$. 
Thus, $S' = S'_1 \cup \dots \cup S'_n$.

By Lemma 5.5, any parallelity handle of
${\cal H}$ that lies in $M'$ is a parallelity handle for
${\cal H}'$. Note that, because each $K_i$ is prime, any incompressible
annulus properly embedded in $M'$ with boundary in $S'$ is boundary parallel.
Also, $\partial M'$ is incompressible, since each $K_i$ is non-trivial.
Hence, by Proposition 5.8, 
${\cal H}'$ has a generalised parallelity bundle ${\cal B}$ 
that contains all parallelity handles of ${\cal H}'$ and that 
consists of $I$-bundles over discs.

We claim that all of $R'$ was removed when constructing ${\cal H}'$
from ${\cal H}$. For, in each annular simplification, each
component of $R$ is either entirely removed or left untouched,
by condition (iv) in Definition 5.4 (the definition of an annular simplification).
Hence, each component of $R$ is either a subset of ${\cal B}$ or
was removed. But each component of $R'$ contains
a component of $\partial M \cap \partial X$ and so cannot lie in an
$I$-bundle over a disc. Thus, each component of $R'$ was
removed in the annular simplifications, as claimed.
Recall that, in this situation, we declared that ${\rm cl}(\partial M' - S')$ is a collection of vertical
boundary components of generalised parallelity bundles that
are removed. In particular, ${\rm cl}(\partial M' - S')$ inherits the
structure of an $I$-bundle. Now, the vertical boundary of ${\cal B}$
lies in a collection of parallelity handles of ${\cal H}'$,
by (iv) in Definition 5.2 (the definition of a generalised parallelity bundle).
None of these handles is incident to $\partial M' - S'$,
for this would have violated condition (iv) in Definition 5.4 (the definition of an annular
simplification). So, the vertical boundary of ${\cal B}$ is disjoint from $\partial M' - S'$.

\vskip 6pt
\noindent {\caps 6.4. Choosing the curves $C_1, \dots, C_n$}
\vskip 6pt

Recall that we are going to pick a simple closed curve $C_i$
on the boundary of each $X'_i$. The union of these curves will
be $K_1 \sqcup \dots \sqcup K_n$, and their projection will be
the diagram $D'$. We first declare that $C_i \cap {\rm cl}(\partial M' - S')$
is a fibre in the $I$-bundle structure on ${\rm cl}(\partial M' - S')$.
We may also arrange that this arc lies in a 0-handle of ${\cal H}'$.
The remainder of $C_i$ will be an arc $\alpha_i$ in $S_i'$ joining
the endpoints of this fibre. We will choose the
arcs $\alpha_i$ so as to control the crossing number of the
resulting diagram $D'$.

Let $\alpha = \alpha_1 \cup \dots \cup \alpha_n$. Our first task is to ensure
that $\alpha$ sits nicely with respect to a certain handle structure
on $S'$, which we now define. Now, $S'$ consists of two parts:
\item{(i)} copies of normal discs of $A$, which we denote by $S'_A$;
\item{(ii)} vertical boundary components of generalised parallelity bundles that have been removed
by annular simplifications, which we denote by $S'_V$.

\noindent Thus, $S' = S'_A \cup S'_V$. The normal discs in $A$ specify a handle structure 
on $S'_A$, where (for $j = 0,1,2$) the $j$-handles are the intersection with
${\cal H}^j$. We may extend this to a handle structure on all of $S'$, by
declaring that the intersection of $S'_V$ with ${\cal H}^0$
is 1-handles (running between the two boundary components of
the relevant vertical annulus), and that the remainder of $S'_V$
is 2-handles. 

We may ambient isotope
$\alpha$ within $S'$, keeping $\partial \alpha$ fixed, so that it misses the 2-handles of $S'$ and 
respects the product structure on the 1-handles. This implies that in a 1-handle $H_1 = D^1 \times D^2$
of the diagrammatic handle structure, $\alpha \cap H_1$ is of the form $D^1 \times P$,
where $P$ is a disjoint union of points in the interior of $D^2$.
We may also arrange that the restriction of the diagrammatic projection map
to $\alpha \cap H_1$ is an embedding. (Recall that the 1-handles of the diagrammatic
handle structure are `horizontal' in $S^3$.)

Now, some handles of $S'$ lie in the generalised parallelity bundle ${\cal B}$.
But, crucially, ${\cal B}$ is a collection of $I$-bundles over discs, disjoint from $\partial M' - S'$.
Thus, the intersection $S' \cap {\cal B}$, which is the horizontal
boundary of ${\cal B}$, 
is a collection of discs in the interior of $S'$. We may therefore
pick the arcs $\alpha_i$ so that they avoid the generalised parallelity bundle ${\cal B}$,
without changing the choice of $\partial \alpha_i$.

Define the {\sl length} of $\alpha$ to be the number of 0-handles
of $S'$ that it runs through (with multiplicity). We pick $\alpha$ so that it has shortest
possible length among arcs that avoid the generalised parallelity bundle ${\cal B}$. 
Hence, for each 0-handle $D$ of $S'$, $\alpha \cap D$ is at most one arc.
Otherwise, we may find an embedded arc $\beta$ in
$D$ such that $\alpha \cap \beta = \partial \beta$ and such
that the endpoints of $\beta$ lie in distinct components of $D \cap \alpha$.
Cut $\alpha$ along $\partial \beta$, discard the arc that misses $\partial S'$ and replace it
by $\beta$. The result is a shorter collection of arcs than $\alpha$, which is a contradiction.

The fact that $\alpha$ misses the generalised parallelity bundle ${\cal B}$ implies, in particular,
that it misses the parallelity regions in $X'_1 \cup \dots \cup X'_n$ between parallel normal
discs of $A$. Thus, in any handle of the diagrammatic handle structure, 
$\alpha$ can run over at most two normal discs of
any given disc type. For, in any collection of three or more parallel normal
discs, all but the outer two discs have parallelity regions
on both sides. Suppose that $\alpha$ runs over a disc $D$ that is not
an outer disc. Then $D$ lies in $X'_1 \cup \dots \cup X'_n$,
because $\alpha$ does. Moreover, one of the parallelity regions adjacent to
$D$ is a parallelity handle of ${\cal H}'$ and so lies in
${\cal B}$. But $\alpha$ misses ${\cal B}$, which is a contradiction.

We now wish to find an upper bound on the crossing number of the
diagram $D'$. In order to do this, we need to be precise about
how the normal discs of $A$ lie in each 0-handle, and how $\alpha$
intersects these discs.

\vskip 6pt
\noindent {\caps 6.5. The position of the normal discs of $A$}
\vskip 6pt

The normal discs of $A$ come in two types: those that miss the boundary
of $X$, and those that intersect $\partial X$. The discs of the latter type
have been removed when creating ${\cal H'}$. Thus, the normal discs
of $A$ that lie in $S'$ all miss $\partial X$, and are therefore
triangles and squares in unexceptional 0-handles, as shown in Figure 11.
Each unexceptional 0-handle of the diagrammatic handle structure simultaneously supports
at most four triangle types and at most one type of square.

We say that a normal disc $E$ properly embedded in a 0-handle of the diagrammatic
handle structure is {\sl flat}
if, for each point $x$ in the diagram 2-sphere, the inverse image of
$x$ in $E$ under the projection map is either empty, a single point, or an
arc in $\partial E$. We say that $E$ is {\sl semi-flat}
if it contains a properly embedded arc $\delta$ such that
the closure of each component of $E - \delta$ is flat.
We say that a flat disc is {\sl convex} if the image of its projection
is a convex subset of the diagram 2-sphere. (Recall that we have assigned a
Euclidean metric to $S^2 - \{ \infty \}$, where $S^2$ is the diagram
2-sphere and $\infty$ is a point that is distant from the crossings,
and hence the 0-handles.) We say that a semi-flat
disc is {\sl piecewise-convex} if its two flat subdiscs are convex.
It is clear that we can make all of the triangles and squares of $A$ 
simultaneously flat and convex, apart from 
certain squares, as shown in Figure 19, which can be made semi-flat
and piecewise convex.

\vskip 18pt
\centerline{
\epsfxsize=3.5in
\epsfbox{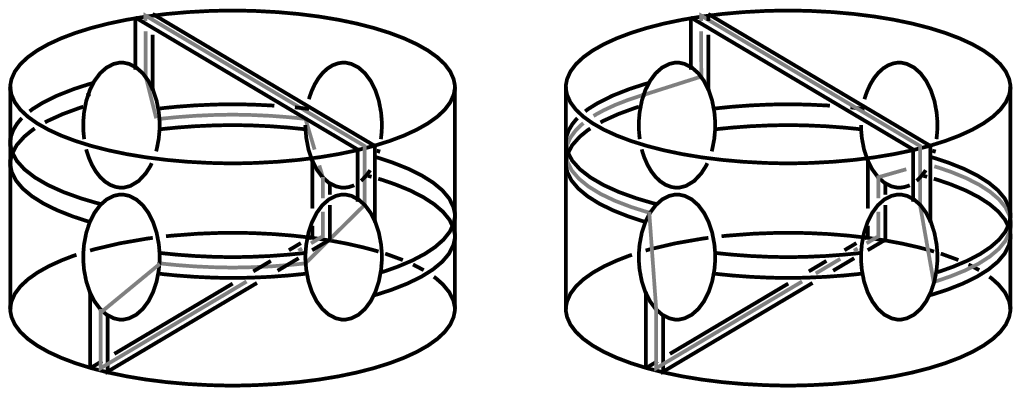}
}
\vskip 6pt
\centerline{Figure 19.}

When a triangle or square $E$ is flat and convex, realise $\alpha \cap E$ as the unique arc in $E$ that projects
to a Euclidean geodesic under the diagrammatic projection map and that
has the same endpoints. When $E$ is semi-flat and piecewise convex, it contains a properly
embedded arc $\delta$ such that each component of $E - \delta$ is
flat and convex. We may therefore realise the intersection of $\alpha$ with
each such component as an arc that projects to Euclidean geodesic.
So, $\alpha \cap E$ projects in this case to the concatenation of
two Euclidean geodesics. We call such an arc {\sl bent}.

Thus, we have defined the precise location of the simple closed curves $C$. Its projection is the
diagram $D'$. We now wish to bound the crossing number of $D'$.

\vskip 6pt
\noindent {\caps 6.6. Justifying the constant $152$}
\vskip 6pt

The arcs $C \cap (\partial M' - S')$ are the interiors of fibres in the parallelity
regions between adjacent normal discs of $A$. Thus, we may clearly
arrange that their projections in $D'$ have disjoint images and
are disjoint from the image of $C \cap S'$.
The crossings of $D'$ therefore arise when arcs of intersection between $C$
and distinct normal discs of $A$ do not have disjoint projections.
The crossings occur only within the projections of the unexceptional
0-handles. Consider one such 0-handle $H_0$. The intersection $C \cap H_0$
lies in at most 10 normal discs, of which all but at most 2 are flat and convex. 
The non-flat discs are semi-flat and piecewise convex. Its intersection
with each flat disc projects to straight arc. The intersection with
each semi-flat disc projects to a bent arc. Now, the projection of
two straight arcs has at most one crossing. There are at most $\left ( 8 \atop 2 \right ) = 28$
such crossings. The projection of two
bent arcs has at most 4 crossings. The projection of a bent arc
and a straight arc has at most 2 crossings, and there are therefore
at most 32 such crossings. So, the number of
crossings of $\alpha$ in the projection of $H_0$ is at most
$28 + 4 + 32 = 64.$

In fact, we may improve this estimate a little. For each type of triangle in $H_0$,
there is another triangle type in $H_0$, with the property that the projections
of these triangles to the diagram are disjoint. Hence, we may reduce the
upper bound on the number of straight-straight crossings by $8$ to $20$.
We also note that the projection of the arc $\delta$ in each semi-flat disc can
each be arranged to lie in the `centre' of the projection of $H_0$.
More specifically, we can define the central region to be
the intersection of the image of the two horizontal pieces of ${\cal F}^1 \cap H_0$,
and we can ensure that $\delta$ is a vertical arc that projects to
this central region. Thus, the two geodesics in each bent arc
are nearly radial. This implies that the projection of each such geodesic
is disjoint from the projection of two triangle types. Hence,
the upper bound on the number of straight-bent crossings can be reduced by $16$ to $16$.
Finally, by carefully arranging the semi-flat discs, we may
ensure that if there are two bent arcs, then they intersect at most twice. So, the
number of crossings of $C$ in the projection of $H_0$ is at most
$20 + 2 + 16 = 38$. An example of the projection of $H_0$ containing
38 crossings is given in Figure 20.

The number of unexceptional 0-handles is $4\ c(D)$. Thus, $c(D')$
is at most $152 \ c(D)$, which proves Theorem 1.1.
$\square$

\vskip 18pt
\centerline{
\epsfxsize=2.3in
\epsfbox{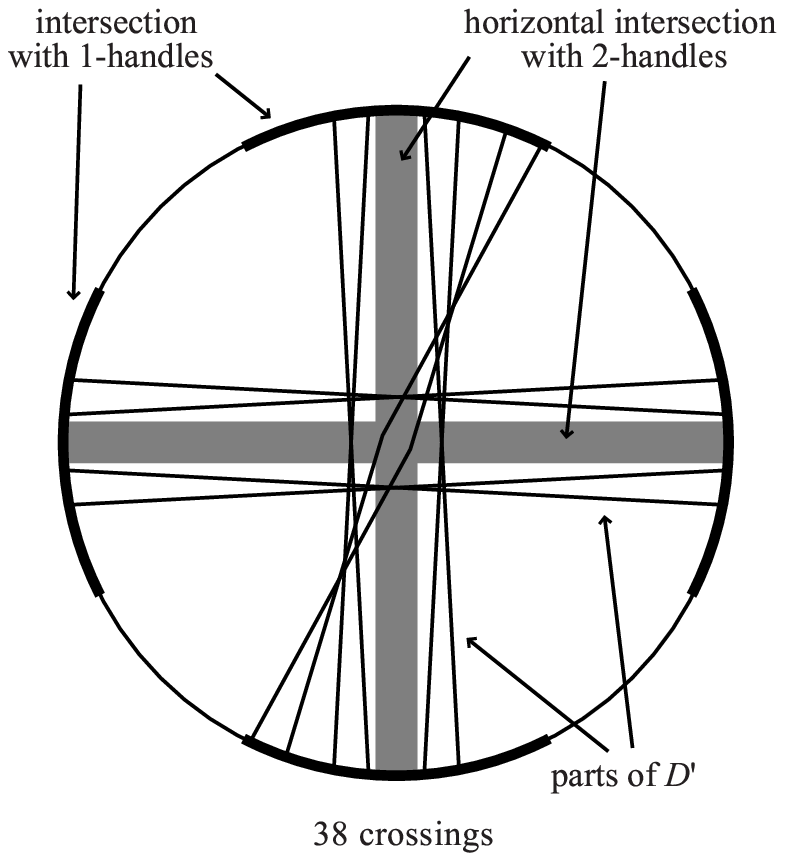}
}
\vskip 6pt
\centerline{Figure 20.}

\vskip 18pt
\centerline{\caps References}
\vskip 6pt

\item{1.} {\caps Y. Diao}, {\sl The additivity of crossing numbers,}
J. Knot Theory Ramifications 13 (2004) 857--866. 

\item{2.} {\caps J. Hass, J. Lagarias}, {\sl The number of Reidemeister moves needed
for unknotting}, J. Amer. Math. Soc. 14 (2001) 399--428.

\vfill\eject
\item{3.} {\caps L. Kauffman}, {\sl State models and the Jones polynomial}, Topology 26 (1987) 395--407.

\item{4.} {\caps M. Lackenby}, {\sl The crossing number of satellite knots}, To appear.

\item{5.} {\caps S. Matveev,} {\sl Algorithmic topology and classification of $3$-manifolds},
Algorithms and Computation in Mathematics, Volume 9, Springer (2003). 

\item{6.} {\caps K. Murasugi}, {\sl The Jones polynomial and classical conjectures in knot theory},
Topology 26 (1987) 187--194.

\item{7.} {\caps M. Thistlethwaite,} {\sl A spanning tree expansion of the Jones polynomial},
Topology 26 (1987) 297--309. 

\vskip 12pt
\+ Mathematical Institute, University of Oxford, \cr
\+ 24-29 St Giles', Oxford OX1 3LB, United Kingdom. \cr

\end